\baselineskip=15pt plus 2pt
\magnification =1200
\nopagenumbers
\voffset=2\baselineskip
\headline={\ifnum\pageno=1 \hss {\smalletters softhardpaper.tex} \hss \else\tt G.Blower: 
Operators for soft and hard edges\ \folio\fi}
\def\sqr#1#2{{\vcenter{\vbox{\hrule height.#2pt\hbox{\vrule width.#2pt
height#1pt\kern#1pt \vrule width.#2pt}\hrule height.#2pt}}}}
\def\square{\mathchoice\sqr64\sqr64\sqr{2.1}3\sqr{1.5}3}
 at 10truept
\font\smalletters=cmr8 at 10truept
\font\smallbold =cmbx10 at 10truept at 12truept

\font\medbold =cmbx10 scaled\magstep2
\centerline {\medbold Operators Associated with Soft and Hard}\par
\centerline {\medbold Spectral Edges from Unitary Ensembles}\par
\vskip.05in
\centerline {\bf Gordon Blower}\par
\vskip.05in
\centerline  {\sl Department of Mathematics and Statistics, 
Lancaster University}\par
\centerline  {\sl Lancaster, LA1 4YF, England, UK. E-mail: 
g.blower@lancaster.ac.uk}\par
\vskip.1in
\centerline {27th April 2006}\par
\vskip.1in
\hrule
\vskip.1in
{{\noindent {\smallbold Abstract.} Using Hankel operators and shift-invariant
subspaces on Hilbert space, this paper develops
the theory of the operators associated with soft and hard edges of
eigenvalue distributions of random matrices.
Tracy
and Widom introduced a projection operator $W$ to describe the
soft edge of the spectrum of the Gaussian unitary ensemble. 
The subspace $WL^2$ is simply invariant under the translation semigroup $e^{itD}$ $(t\geq 0)$
and invariant under the Schr\"odinger semigroup $e^{it(D^2+x)}$ $(t\geq 0)$;  
these properties characterize $WL^2$ via Beurling's theorem. The Jacobi ensemble
of random matrices has positive eigenvalues which tend to accumulate near to the hard edge
at zero. This paper identifies a pair of unitary groups that satisfy the von
Neumann--Weyl anti-commutation relations and leave invariant certain subspaces of
$L^2(0,\infty )$ which are invariant for operators with Jacobi
kernels. Such Tracy--Widom operators are reproducing kernels for
weighted Hardy spaces, known as Sonine spaces. Periodic
solutions of Hill's equation give a new family of Tracy--Widom type operators.\par}} 
\vskip.05in

\par
\vskip.05in

\par
\vskip.1in
\hrule
\vskip.1in

\noindent {\bf 1 Introduction}\par
\indent This paper concerns the spectral theory and invariant
subspaces of operators that arise in random matrix theory,
particularly the soft and hard edges that occur on the limiting
eigenvalue distributions of the Gaussian and Jacobi unitary
ensembles. Tracy and Widom [28, 29, 30] introduced various operators to describe
the soft edge of the spectrum of the Gaussian unitary ensemble; that
is, the eigenvalues near to the supremum of the support of the 
equilibrium distribution. Burnol proposed that the theory of random matrices should be
expressed in terms of Sonine spaces [7, p 692].
Here we develop this theory in a systematic
manner to show that the Tracy--Widom calculations are instances of
more general
results on Hankel operators, and introduce new settings where the theory applies. \par
\indent In section 2 we consider operators on $L^2({\bf R})$ with kernels
$$W(x,y)={{A(x)B(y)-B(x)A(y)}\over {x-y}}\eqno(1.1)$$
\noindent where ${\hbox{col}}\, [A,B]$ satisfies a first-order linear differential
equation, and  give sufficient conditions for $W$ to be the square of a Hankel 
integral operator.
Further, we show that the determinants $\det (I-zW)$ are related to the solutions of Marchenko
integral equations. As we show in section 3, kernels such as $W$ arise as reproducing
kernels for weighted Hardy spaces on the upper half-plane ${\bf C}_+=\{ z: \Im
z>0\}$ as in [2, 6].\par
Let $J$ be the `flip' map
$Jf(x)=f(-x),$ and $M_u$ the multiplication operator $f\mapsto uf$. The
classical Hardy space $H^2$ consists of the holomorphic functions $F$
on ${\bf C}_+$ such that $\sup_{y>0}\int_{-\infty}^\infty
\vert F(x+iy)\vert^2\, dx<\infty$, and we identify such a function
with its $L^2$ boundary values. The Fourier transform is ${\cal F}f(\xi
)=\int_{-\infty}^\infty e^{-ix\xi} f(x)\, dx/\sqrt{2\pi}$. Given $u\in L^\infty$,
the bounded linear operator $\sqrt{2\pi} {\cal F}^*M_u{\cal F}^*$ is the Hankel
operator $\Gamma_{u}$ on $L^2({\bf R}_+)$ that has distributional kernel ${\cal F}^*u(x+y)$ as in
[23]. \par

\indent We recall how examples of such operators appear in the theory of the Gaussian unitary
ensemble. Let $x_{j,k}$ and $y_{j,k}$ $(1\leq j\leq k\leq n)$ be a family of
mutually independent $N(0,1/n)$ random variables. We let $X_n$ be the
$n\times n$ Hermitian matrix that has entries
$[X_n]_{jk}=(x_{j,k}+iy_{j,k})/\sqrt 2$ for $j<k$,
$[X_n]_{jj}=x_{jj}$ for $1\leq j\leq n$ and 
$[X_n]_{kj}=(x_{j,k}-iy_{j,k})/\sqrt 2$ for $j<k$; the space of all such
matrices with the probability measure $\sigma_n^{(2)}$ forms the 
{\sl Gaussian unitary ensemble}.\par
\indent {\sl Bulk of the spectrum.} The eigenvalues of $X_n$
are real and may be ordered as $\lambda_1\leq \dots \leq \lambda_n$,
so their positions are specified  by the empirical
distribution $\mu_n =(1/n)\sum_{j=1}^n \delta_{\lambda_j}$. As 
$n\rightarrow\infty$, the empirical distributions converge weakly
to the equilibrium distribution, namely the Wigner semicircle law 
$$\rho (dx)={{1}\over {2\pi}}{\bf I}_{[-2,2]}(x)\sqrt {4-x^2} dx,\eqno(1.2)$$
\noindent for almost all
sequences $(X_n)$ of matrices under $\otimes_{n=1}^\infty \sigma_n^{(2)}$. The bulk of the spectrum
consists of those eigenvalues in $[-2,2].$ See [22, p 93]. \par
\indent Let $B_t$ 
be the operator on $L^2({\bf R})$ that has kernel
$$B_t(x,y)={{\sin t\pi (x-y)}\over {\pi (x-y)}},\eqno(1.3)$$
let ${\bf I}_S$ be the indicator function of a set $S$, and let $P_{(\alpha ,\beta )}$ 
be the orthogonal projection on
$L^2({\bf R})$ given by $P_{(\alpha ,\beta )}f(x)={\bf I}_{(\alpha
,\beta )}(x)f(x)$; for brevity we write $P_+=P_{(0,\infty )}$ and
$P_-=P_{(-\infty ,0)}.$\par 
\indent Let $E_{\sigma_n}(k;\alpha ,\beta )$ be
the probability with respect to $\sigma_n^{(2)}$ that $(\alpha ,
\beta )$ includes exactly $k$ eigenvalues. Mehta and Gaudin [22] showed that  
$$E_{\sigma_n}(k;\alpha , \beta )\rightarrow 
{{(-1)^k}\over {k!}}\Bigl( {{d^k}\over {dz^k}}\Bigr)_{z=1}
\det \Bigl[ I-zP_{(\alpha , \beta )}B_1
P_{(\alpha , \beta )}\Bigr].\eqno(1.4)$$ 
\indent This determinant can alternatively be expressed in terms of the operator\par
\noindent $\Psi_a:L^2[-a,a]\rightarrow L^2$ that has 
kernel $\Psi_a(x,y)=e^{ixy}{\bf I}_{[-a,a]}(y)/\sqrt {2\pi}$ and satisfies $\Psi_a\Psi_a^*=B_{a/\pi}$.\par
\vskip.05in
\indent {\sl Hard edges.} Let $Y_n$ random $n\times n$ matrices with independent $N(0,1/n)$ entries, 
and let $0\leq \lambda_1\leq \dots
\leq\lambda_n$ be the eigenvalues of the positive operator $Y_n^*Y_n$. Then
$\nu_n={{1}\over {n}}\sum_{j=1}^n \delta_{\lambda_j}$ is the empirical eigenvalue
distribution, and $\nu_n$ converges weakly almost surely to the
Marchenko--Pastur distribution; so that
$$\int_0^\infty f(x)\nu_n(dx)\rightarrow \int_0^\infty f(x)\sqrt{
{{4-x}\over {x}} }\, {{dx}\over {2\pi}}\eqno(1.5)$$
\noindent almost surely as $n\rightarrow\infty$ for all continuous and bounded
real functions $f$. Thus the $\lambda_j$ tend
to accumulate near to their minimum possible value of zero, where the density of the
limiting distribution is unbounded; this is the hard edge effect.\par
\indent Hard edges also arise from random matrices of the Jacobi and Laguerre
ensembles. The Jacobi ensemble of order $N$ with parameters $\nu ,\gamma >-1/2$ at inverse
temperature $\beta >0$ is the joint
distribution function
$$\sigma^{(\beta )}_{N,J} (dx)={{1}\over {Z_N}}\prod_{j=1}^N (1+x_j)^{\beta \gamma}
(1-x_j)^{\beta \nu }\prod_{1\leq j<k\leq N} (x_k-x_j)^\beta\, dx_1\dots
dx_N\eqno(1.6)$$
\noindent where $-1\leq x_1\leq \dots \leq x_N\leq 1$ are the eigenvalues.\par 
\indent Forrester [12] showed that the integral operator $F^{a,b}$ on $L^2((0,1),dx)$ with 
kernel
$$F^{a,b}(x,y)={\bf I}_{(a,b)}(x){{J_\nu (\sqrt x)\sqrt{y}J_\nu' (\sqrt y)-\sqrt{x}J_\nu' (\sqrt x)
J_\nu (\sqrt y)}\over {2(x-y)}}{\bf I}_{(a,b)}(y)\eqno(1.7)$$
\noindent determines the limiting distribution of scaled eigenvalues $x_j/(4n)$ from the Laguerre
ensemble near to the hard edge, and conjectured that a similar result holds for the
Jacobi ensemble. Using the orthogonal polynomial technique, Forrester and Rains [13] have verified the cases of $\beta =1,2$ and $4$
following earlier work by Borodin [5] and Due\~nez.\par 
\indent We introduce the scaled eigenvalues $\xi_j$ by $x_j=\cos \xi_j/\sqrt N$,
to ensure that the mean spacing of the $\xi_j$ is of order $O(1)$ near to the hard edge 
at $x_j\approx 1$. One can show that 
$$\sigma_{N,J}^{(2)} [ (a,b)\quad{\hbox{contains no}}\quad \xi_j]\rightarrow \det
(I-F^{a,b})\qquad (N\rightarrow  \infty ).\eqno(1.8)$$  
For subsequent analysis we change variables by writing $x=e^{-2\xi}$ and $y=e^{-2\eta }$ so
that $\xi ,\eta \in (0,\infty )$ for $x,y\in (0,1)$. 
 Let $G_\ell$ be the unitary integral operator on $L^2({\bf R})$ that has
kernel $e^{-\ell -\xi-\eta }J_\nu (e^{-\ell -\xi-\eta });$ let $Q_\ell
=G_\ell P_+ G_\ell$ $(\ell \in {\bf R})$, which gives a strongly continuous
family of orthogonal projections. For compact operators $S$ and $T$ on Hilbert space, the
spectrum of $ST$ equals the spectrum of $TS$.
The integral operator $\Phi_\ell =P_+G_\ell P_+$ on
$L^2(0,\infty )$  is
Hilbert--Schmidt, and when $0<a<1$ and $\alpha =-(1/2)\log a$ satisfies
$$\det (I-zF^{0,a})=\det (I-z\Phi_{(\alpha )}^2).\eqno(1.9)$$
\indent In section 4 we interpret these operators on the Sonine spaces
$u_\nu H^2$ where $u_\nu (x)=2^{ix}\Gamma((1+\nu +ix)/2)/\Gamma ((1+\nu -ix)/2)$.\par
\vskip.05in
\indent {\sl Soft edge of the spectrum.} We recall some results of Tracy and Widom
[12, {28}] concerning the largest few
eigenvalues. The Airy function
 ${\hbox{Ai}}(x)$, as defined by the oscillatory integral
$${\hbox{Ai}}(z)={{1}\over{2\pi}}\int_{-\infty}^\infty
e^{i(zt+t^3/3)}\, dt,\eqno(1.10)$$
\noindent satisfies 
the Airy differential equation [27, page 18] $y''-xy=0.$ Let $W_{1/3}$ be the integral operator
on $L^2({\bf R})$ defined by the {\sl Airy kernel}
$$W_{1/3}(x,y)=
{{{\hbox{Ai}}(x){\hbox{Ai}}'(y)-{\hbox{Ai}}'(x){\hbox{Ai}}(y)}\over
{x-y}}.\eqno(1.11)$$
\noindent We scale the eigenvalues of the Gaussian ensemble by introducing
$$\xi_j=n^{2/3}\Bigl( {{\sqrt {2}}\over {\sqrt{n}}}\lambda_j-2\Bigr),
\eqno(1.12)$$
\noindent and let $E_{\sigma_n}(k;\xi ; \alpha , \beta )$ be the probability
with respect to $\sigma_n^{(2)}$ that $(\alpha ,\beta )$ contains exactly
$k$ of the $\xi_j$ $(j=1, \dots , n)$; see
[22, page 116, A7]. Aubrun [3] proved that the operator 
$W_{1/3}^{\alpha, \beta }=P_{(\alpha ,\beta )} W_{1/3}P_{(\alpha
,\beta )}$ on $L^2({\bf R}_+)$  is of
trace class for $0<\alpha <\beta \leq\infty$, and 
$$E_{\sigma_n}(k; \xi ;\alpha ,\beta )\rightarrow {{(-1)^k}\over {k!}}\Bigl({{d^k}\over
{dz^k}}\Bigr)_{z=1} \det\bigl(
I-zW_{1/3}^{\alpha ,\beta}\bigr)\qquad
 (n\rightarrow\infty ).\eqno(1.13)$$
\indent The compression of $W_{1/3}^{\alpha , \infty}$ to $L^2(\alpha ,\infty )$ may
be identified, under the change of variables $s\mapsto \alpha +s,$ with 
$\Gamma_{(\alpha )}^2$ where the Hankel integral operator $\Gamma_{(\alpha )}$ on 
$L^2[0,\infty )$ satisfies
$$\Gamma_{(\alpha)} f(s)=\int_0^\infty {\hbox{Ai}}(\alpha +s+t)f(t)\,dt\qquad (f\in L^2(0,\infty
)).\eqno(1.14)$$
\noindent The spectrum of $P_{(\alpha ,\beta
)}\Gamma_{(0)}^2P_{(\alpha ,\beta )}$ equals the spectrum of
$\Gamma_{(0)}P_{(\alpha ,\beta )}\Gamma_{(0)}$, and hence 
$$\det (I-zP_{(\alpha ,\infty )}W_{1/3}P_{(\alpha ,\infty )})=
\det( I-z\Gamma_{(\alpha )}^2).\eqno(1.15)$$
\vskip.05in
\indent {\sl Edge distributions and KdV.} For $0\leq t\leq 1$ let $w(x;t)$ be the unique solution to the Painlev\'e II
equation $w''=2w^3+xw$ that satisfies $w(x;t)\asymp -\sqrt {t}{\hbox{Ai}}(x)$ as
$x\rightarrow\infty$. By the theory of inverse scattering for the concentric
Korteweg--de Vries equation, this solution is given by the Fredholm determinant
$$w(x;t)^2=-{{\partial^2}\over {\partial x^2}}\log\det( I-t\Gamma_{(x)}^2);\eqno(1.16)$$
\noindent see [10, p. 86, 174]. 
Tracy and Widom [{28}] introduced the cumulative
distribution function $F(x;t)=\det (I-t\Gamma_{(x)}^2)$ so that 
$$F(x;t)=\exp\Bigl( -\int_x^\infty (y-x)w(y;t)^2\, dy\Bigr);\eqno(1.17)$$
\noindent in particular, $F(x;1)$ is the {\sl Tracy--Widom distribution.} \par
\indent {\sl Operators for parts of the spectrum.} In section 5 we show
how $W_{1/3}$ and $\Gamma_{(0)}$ arise from the Airy group $e^{itD^3}$ on
$L^2({\bf R})$ where
$D=-i{{\partial}\over {\partial x}}$. 
\noindent By suitable changes of variable we arrange that the edge of
the support of the equilibrium distribution is at zero and we consider
the operators on $L^2(\alpha ,\infty )$ that describe the probability
that scaled eigenvalues lie in $(\alpha ,\infty )$.\par 
\indent  The relative positions of $L^2({\bf R}_+)$ and $W_{1/3}L^2({\bf
R})$ are described by $P_+W_{1/3}P_+$. Generally we have
$W_{t}=e^{itD^3}P_-e^{-itD^3}$ and the complementary 
orthogonal projection is $W_{t}^\perp =I-W_{t}=e^{itD^3}P_+e^{-itD^3}$. For 
comparison, $R_+={\cal F}^*P_+{\cal F}$ and $R_-={\cal F}^*P_-{\cal F}$ are the
Riesz projections on $L^2$ that have images $H^2$ and $\overline{H^2}$ respectively. These formul\ae\quad suggests an analogue of prediction theory such that the subspace $W_{t}^\perp L^2({\bf R})$ 
corresponds to the subspace $L^2({\bf R}_+)$ which represents the future,
and such that the unitary operator $e^{itD^3}$ plays the rol\^e corresponding to
the inverse Fourier transform ${\cal F}^*$ in the Hardy space theory
of [17]. The projections $\tau_\alpha
W_t^\perp\tau_{-\alpha}=e^{itD^3}P_{(\alpha ,\infty )}e^{-itD^3}$ form
a decreasing nest as $\alpha$ increases. The spectrum of the Hankel operator determines the
limiting eigenvalue distribution via (1.13) and (1.15).\par
\indent The following table describes the analogy between the operators and
subspaces in the various cases.\par
$$\matrix{ {}& {\hbox{Classical}}&{\hbox{Bulk}}&{\hbox{Hard Edge}}&
{\hbox{Soft Edge}}\cr
{\hbox{Future Projection}}& {\cal F}^*P_+{\cal F}&{\cal F}^*P_{(-a,a)}{\cal F}&
G_\ell P_+G_\ell &e^{iD^3/3}P_+e^{-iD^3/3}\cr
{\hbox{Future space}}& H^2& B_{a/\pi}L^2&Q_\ell L^2& W_{1/3}^\perp L^2\cr
{\hbox{Subspace position}}&{}&e^{i2ax}H^2\subset H^2& u_\nu \overline{H^2}\cap H^2\neq 0&
e^{itx^3}H^2\cap H^2=0\cr
 {\hbox{Painlev\'e Equation}}& {}&\sigma {\hbox{-P}}_{{\hbox{\smalletters V}}}&
{\hbox{P}}_{{\hbox{\smalletters III}}} &{\hbox{P}}_{{\hbox{\smalletters II}}}\cr
{\hbox{Hankel operator}}&{}&\Psi_a&\Phi_\ell
&\Gamma_{(0)}\cr}$$
\vskip.05in
\indent {\sl Weyl Relations and Invariant Subspaces.}\par
\vskip.05in
\noindent {\bf Definition.} Let $(U_t)$ $(t\in {\bf R})$ be a $C_0$
 (strongly continuous)
group of unitary operators on an infinite-dimensional separable Hilbert space ${
H}$, and let ${K}$ be a closed linear subspace of $H$. We say that $K$ is {\sl doubly
invariant} for $(U_t)$ when $U_tK\subseteq K$ for all $t\in {\bf R}$. 
Further, ${K}$ is 
{\sl simply invariant} for $(U_t)$ $(t\geq 0)$ when $U_tK\subseteq K$ for
$t\geq 0$ and moreover $\cap_{t\geq 0}U_t{K}=\{0\}$.\par 
\vskip.05in
\indent Beurling and Lax characterized the subspaces of $L^2$ that are
invariant for the shift operators $S_s: f(x)\mapsto e^{isx} f(x)$; see [15, 17 page 114
]. For notational simplicity, we sometimes write $e^{i\ell x}H^2=\{ e^{i\ell x}
f(x): f\in H^2\}$.\par
 A closed linear subspace ${\cal T}$ is simply invariant for the
semigroup $\{ S_s :s \geq 0\}$, if and only if there exists a
unimodular measurable function $u$ such that ${\cal T}=uH^2=\{ uf:f\in
H^2\}$; such a $u$ is uniquely determined up to a unimodular constant
factor. In each case we start by making unitary transformations to identify $u$
and to determine the relative positions of $uH^2$ and $H^2$. 
Either $uH^2\cap H^2=0$ or there exist inner functions $v$ and $w$,
uniquely determined up to a unimodular constant factor, such that $u=v\bar w$,
$uH^2\cap H^2=vH^2$ and $vwH^2=vH^2\cap wH^2$. In sections 4 and 5, we find $uH^2\cap
H^2=0$, so we factorize $u(z)=\overline{E(\bar z)}/E(z)$ where $E$ is a
meromorphic function on ${\bf C}$ that has no zeros. By using de Branges's version
of Beurling's theory [6], we are able to show that
$W$ from (1.1) is unitarily equivalent to $\Gamma^*_{u^*}\Gamma_{u^*}$ 
and hence that $W$ is the 
reproducing kernel of some  
weighted Hardy spaces of holomorphic functions inside ${\bf C}_+$.\par
\vskip.05in
\noindent {\bf Definition.} A Weyl pair $(U_s,V_t)$ consists of a pair of $C_0$
unitary groups $(U_s)$ $(s\in {\bf R})$ and  $(V_t)$ $(t\in {\bf R})$ on $H$ that satisfy $U_sV_t=e^{ist} V_tU_s$ for all $(s,t\in {\bf R}).$\par
\vskip.05in
\indent The shifts 
$S_s$ $(s\in {\bf R})$ and the translations $\tau_t=e^{-itD}$ $(t\in {\bf R})$ give a Weyl pair on $L^2$; moreover, this is the unique
representation of the Weyl relations of multiplicity one on $L^2$, up to unitary 
equivalence; see [32]. Katavolos and Power [16] obtained the following description of the 
invariant subspaces for a Weyl pair of multiplicity one. Let ${\cal L}$ be the space of 
orthogonal projections $P$ onto closed linear subspaces $K$
of $L^2$ that are invariant under $(S_t)$ $(t\geq 0)$ and $(\tau_{-t})$ $(t\geq 0)$, 
where ${\cal L}$ has the strong operator topology. There is a
homeomorphism $\rho :\{ z: \vert z\vert \leq 1\}\rightarrow {\cal L}$ such that 
:\par
\indent (1) $\rho (-i)=0,$ and $\rho (i)=I$;\par
\indent (2) $\rho (z)L^2$ with $\vert z\vert <1$ is simply invariant 
for both $(S_s)$ $(s\geq 0)$ and $(\tau_{-t})$
$(t\geq 0)$;\par
\indent (3) $\rho (e^{i\theta })L^2$ with $-\pi /2<\theta <\pi /2$ is simply invariant
for $(S_s)$ $(s\geq 0)$ and doubly invariant for $(\tau_{-t})$ $(t\in {\bf R});$\par
\indent (4) $\rho (e^{i\theta })L^2$ with 
$-\pi /2<\pi -\theta <\pi /2$ is doubly invariant for $(S_s)$ $(s\in {\bf R})$ and 
simply invariant
for $(\tau_{-t})$ $(t\geq 0)$.\par
\indent In section 4 we introduce for the Jacobi ensemble an appropriate Weyl pair
for the subspaces $Q_\ell L^2$. On account of the natural ordering of the subspaces 
$Q_\ell L^2$, and the probabilistic interpretation of (1.9), we naturally take the
translations to be one of the groups in the Weyl pair; whereas we need to hunt down
the other one. For the soft edge ensemble, the Weyl pair consists of
the translations $e^{isD}$ and the Schr\"odinger group $e^{it(D^2+x)}$, as we discuss
in section 5.\par
\indent In section 6 we extend these ideas to a new context, namely the
Mathieu functions. Here the KdV equation is $2\pi$ periodic and associated with an
infinite-dimensional manifold. Whereas we do not propose that this
corresponds to a 
natural
random matrix ensemble, the results illustrate the scope of the theory of
Tracy--Widom operators.\par
\vskip.05in
\vskip.05in
\noindent {\bf 2. Kernels from differential equations and the Marchenko integral equation}\par
\indent In this section we prove some results concerning Tracy--Widom operators which
are already known in specific cases from [8, 28, 29, 30]. Here ${\hbox{B}}(H),$ $c^2$ and $c^1$
respectively denote the bounded, Hilbert--Schmidt and trace-class linear operators 
on Hilbert
space $H$.\par
\vskip.05in
\noindent {\bf Lemma 2.1.} {\sl Suppose that $A$ and $B$ are bounded, measurable and real
functions. Then 
$$W(x,y)={{A(x)B(y)-A(y)B(x)}\over {x-y}}\eqno(2.1)$$
\noindent defines a self-adjoint and bounded linear operator on $L^2({\bf R})$.}\par
\vskip.05in
\noindent {\bf Proof.} This follows from the fact that $M_A$, $M_B$ and $R_+$ are
bounded on $L^2$.\par
\rightline {$\square$}\par
\vskip.05in
\noindent {\bf Proposition 2.2.} {\sl Suppose that $A$ and $B$ are bounded, continuous
and integrable functions such that 
$${{d}\over {dx}}\left[\matrix{ A(x)\cr B(x)}\right] =
\left[\matrix{ \alpha (x)&\beta (x)\cr -\gamma (x)& -\alpha (x)}\right]
\left[\matrix{ A(x)\cr B(x)}\right],\eqno(2.2)$$
\noindent where $\alpha ,\beta$ and $\gamma$ are linear functions such that 
$$C={{1}\over{x-y}}\left[\matrix{ \gamma (x)-\gamma (y)&\alpha (x)-\alpha (y)\cr 
\alpha (x) -\alpha (y)& \beta (x)-\beta (y)}\right]
=\left[\matrix{ c&a \cr a&b}\right]\eqno(2.3)$$
\noindent is a negative (semi-) definite constant matrix.\par
\indent (i) Then there exist continuous real
functions $F$ and $G$ such that\/} 
$${{A(x)B(y)-A(y)B(x)}\over {x-y}}=\int_0^\infty \bigl( F(x+t)F(t+y) 
+ G(x+t)G(t+y)\bigr)\, dt\eqno(2.4)$$
\noindent {\sl so $P_+WP_+$ is a sum of two squares of
Hankel operators. \par
\indent (ii) In particular when $C$ has rank one, the operator $P_+WP_+$ on $L^2(0,\infty )$ is
the square of a self-adjoint Hankel operator.}\par
\vskip.05in
\noindent {\bf Proof.} We take a real bilinear pairing and write 
$$ A(x)B(y)-A(y)B(x)=\Bigl\langle \left[\matrix{ 0&-1\cr 1&0}\right]
\left[\matrix{ A(x)\cr B(x)}\right], \left[\matrix{ A(y)\cr
B(y)}\right]\Bigr\rangle,\eqno(2.5)$$
\noindent so by a short calculation
$$\Bigl({{\partial }\over {\partial x}} + {{\partial }\over {\partial y}}\Bigr) 
{{A(x)B(y)-A(y)B(x)}\over {x-y}}= 
\Bigl\langle \left[\matrix{ c&a\cr a&b}\right]
\left[\matrix{ A(x)\cr B(x)}\right], \left[\matrix{ A(y)\cr B(y)}\right]\Bigr\rangle
.\eqno(2.6)$$
\noindent We introduce a real symmetric matrix $X$ such that $X^2=-C$, 
let ${\hbox{col}}\,[\cos\theta ,\sin\theta ]$ and 
${\hbox{col}}\,[-\sin\theta , \cos\theta ]$ be the unit eigenvectors of $X$ corresponding to
eigenvalues $\lambda_1\geq 0$ and $\lambda_2\geq 0$ respectively. Then 
$$\eqalignno{F(x)&= \lambda_1 \bigl( A(x)\cos\theta +B(x)\sin\theta \bigr) ,\cr
G(x)&= \lambda_2 \bigl( -A(x)\sin\theta +B(x)\cos\theta \bigr)&(2.7)}$$
\noindent are bounded, continuous and integrable functions that satisfy 
$$\eqalignno{\Bigl({{\partial }\over {\partial y}} + {{\partial }\over {\partial x}}
\Bigr)& 
{{A(x)B(y)-A(y)B(x)}\over {x-y}}\cr
&= -F(x)F(y)-G(x)G(y)\cr
&=\Bigl({{\partial }\over {\partial x}} + {{\partial }\over {\partial
x}}\Bigr)\int_0^\infty \bigl( F(x+t)F(y+t)+ G(x+t)G(y+t)\bigr)\, dt.&(2.8)}$$
\noindent Hence both sides of (2.4) differ by $f(x-y)$ for some differentiable function
$f$; but both sides converge to zero as $x$ or $y$ tend to $\infty$; so $f=0$,
and equality holds.\par
\rightline {$\square$}\par
\vskip.05in
\indent This result enables us to calculate a determinant as in (1.4),
(1.9) and (1.15).\par
\vskip.05in

\noindent {\bf Theorem 2.3.} {\sl Let } ${\hbox{A}}:(0,\infty )\rightarrow {\bf R}$ {\sl be a continuous 
function such that} \par
\noindent $\int_0^\infty u{\hbox{A}}(u)^2\, du\leq 1.$ {\sl Then}
$$W(u,v)=\int_0^\infty {\hbox{A}}(u+t){\hbox{A}}(t+v)\, dt\eqno(2.9)$$
\noindent {\sl is the kernel of a trace-class operator on $L^2(0,\infty )$ such 
that, when $\vert\kappa\vert <1$,
$$K(x,z)-\kappa^2\int_x^\infty K(x,y)W(y,z)\, dy=\kappa W(x,z)\eqno(2.10)$$
\noindent has a solution $K(x,z)$, which is a trace-class kernel, such that}
$${{\partial }\over {\partial x}}\log\det (I-\kappa^2 P_{(x,\infty )}WP_{(x,\infty )})=\kappa 
K(x,x)\qquad (x>0).\eqno(2.11)$$

\vskip.05in
\indent In the subsequent proof, we shall use the self-adjoint Hankel operator 
$\Gamma_{(x)}$ as in 
$$\Gamma_{(x)}f(t)=\int_0^\infty {\hbox{A}}(x+t+u)f(u)\, du\qquad (f\in L^2(0\,, \infty
)),\eqno(2.12)$$
\noindent where $\Gamma_{(x)}^2$ on $L^2(0,\infty )$ is unitarily
equivalent to $P_{(x,\infty )}WP_{(x,\infty )}$ on $L^2(x, \infty )$.\par
\vskip.05in
\noindent {\bf Lemma 2.4.} {\sl For $\vert
\kappa\vert^2\int_0^\infty u{\hbox{A}}(u)^2\, du<1$ there exists a solution $L\in
L^2((0,\infty )^2)$  to the integral equation}
$$L(x,s)-\kappa^2 \int_0^\infty L(x,y)\int_x^\infty {\hbox{A}}(y+u){\hbox{A}}
(u+s)\, dudy =\kappa {\hbox{A}}(x+s).\eqno(2.13)$$
\vskip.05in
\noindent {\bf Proof of Lemma 2.4.} By the Hilbert--Schmidt theorem applied to
$\Gamma_{(0)}$, there exist $(\varphi_j)$, 
an orthonormal basis of
$L^2(0,\infty )$, and real $\gamma_j$ such that
${\hbox{A}}(x+y)=\sum_{j=1}^\infty \gamma_j
\varphi_j(x)\varphi_j(y)$ and 
$$\sum_{j=1}^\infty
\gamma_j^2=\int_0^\infty\!\!\!\int_0^\infty{\hbox{A}}(u+v)^2\,dudv=\int_0^\infty
u{\hbox{A}}(u)^2\, du<\infty .\eqno(2.14)$$ 
Our solution separates into series $L(x,s)=\sum_{j=1}^\infty
\chi_j(x)\varphi_j(s)$, where ${\hbox{col}}\,[\chi_j(x)]$ satisfies the equation
with column vectors in $\ell^2$
$$ \bigl[ I-\kappa^2\Phi (x)\bigr]{\hbox{col}}\,[\chi_j(x)]=\kappa
\,{\hbox{col}}\,[\gamma_j\varphi_j(x)]\eqno(2.15)$$
\noindent and with matrix
$$\Phi (x)=\Bigl[ \gamma_j\gamma_k\int_x^\infty \varphi_j(y)\varphi_k(y)\,
dy\Bigr]_{1\leq j,k<\infty}.\eqno(2.16)$$
\noindent By the Cauchy--Schwarz inequality, $\Phi (x)$ defines a Hilbert--Schmidt 
operator on $\ell^2$ with norm 
$$\Vert \Phi (x)\Vert_{c^2}^2\leq \sum_{j,k=1}^\infty 
\gamma_j^2\gamma_k^2\Bigl(\int_x^\infty \vert\varphi_j(y)\vert^2\,
dy\Bigr) \Bigl(\int_x^\infty \vert\varphi_k(y)\vert^2\, dy\Bigr)\leq \Bigl( \sum_{j=1}^\infty \gamma_j^2\Bigr)^2;\eqno(2.17)$$
\noindent hence, $I-\kappa^2\Phi (x)$ is invertible whenever $\vert\kappa\vert^2 \Vert
\Phi (x)\Vert_{c^2}<1$. We deduce that (2.15) has a unique solution and by orthogonality 
$$\eqalignno{\int_0^\infty \int_0^\infty \vert L(x,y)\vert^2\, dxds&=\int_0^\infty 
 \sum_{j=1}^\infty \vert \chi_j(x)\vert^2\, dx\cr
&\leq \int_0^\infty \Vert (I-\kappa^2\Phi
(x))^{-1}\Vert^2_{{\hbox{\smalletters B}}(\ell^2)}\Vert 
(\kappa\gamma_j\varphi_j(x))\Vert^2_{\ell^2}\, dx\cr
&\leq \Bigl( 1-\vert\kappa\vert^2\sum_{j=1}^\infty 
\gamma_j^2\Bigr)^{-2}\vert\kappa\vert^2\sum_{j=1}^\infty
\gamma_j^2.&(2.18)\cr}$$
\rightline {$\square$}\par
\noindent {\bf Remark.} When $\sum_{j=1}^\infty \vert\gamma_j\vert <\infty$, the
operator $\Phi (x)$ is trace class and the determinant 
$\det (I-\kappa^2\Phi
(x))$ defines an entire function of $\kappa $.\par
\vskip.05in
\noindent {\bf Proof of Theorem 2.3.} By Lemma 2.4, the integral kernels 
$$W(x,y)=\int_0^\infty {\hbox{A}}(x+s){\hbox{A}}(s+y)\, ds\eqno(2.19)$$
\noindent and
$$K(x,z)=\int_0^\infty L(x,s){\hbox{A}}(s+z)\, ds\eqno(2.20)$$
\noindent are trace class and satisfy the Marchenko integral equation
$$K(x,z)-\kappa^2\int_x^\infty K(x,s)W(s,z)\, ds=\kappa W(x,z).\eqno(2.21)$$ 
\noindent By iterated substitution, we deduce that 
$$K(x,z)=\kappa W(x,z)+\kappa^3W(x, \,.\,) P_{(x,\infty )}(I-\kappa^2
P_{(x,\infty )}WP_{(x,\infty )})^{-1}P_{(x,\infty )}W(\,.\, ,z).$$
\noindent We deduce, by simplifying the integral kernels, that
$$\kappa K(x,x)=\kappa^2\bigl\langle (I-\kappa^2\Gamma^2_{(x)})^{-1}
{\hbox{A}}(x+\,.\, ), {\hbox{A}}(x+\,.\, )\bigr\rangle_{L^2(0,\infty )}.\eqno(2.22)$$ 
\noindent By unitary equivalence as in (2.12), we have
$$\det \bigl(I-\kappa^2 P_{(x,\infty )}WP_{(x,\infty
)}\bigr)=\det \bigl(I-\kappa^2 \Gamma_{(x)}^2\bigr)$$
\noindent and hence 
$${{\partial }\over {\partial x}}\log\det \Bigl( I-\kappa^2 P_{(x,\infty )}WP_{(x,\infty
)}\Bigr)=-\kappa^2{\hbox{trace}}\, \Bigl( (I-\kappa^2 \Gamma_{(x)}^2)^{-1}
{{\partial }\over {\partial x}}\Gamma_{(x)}^2\Bigr),\eqno(2.23)$$
where ${{\partial }\over {\partial x}}\Gamma_{(x)}^2$ is the rank-one
operator that has kernel 
$$\eqalignno{{{\partial }\over {\partial x}}\Gamma_{(x)}^2(u,v)&=
{{\partial }\over {\partial x}}\int_0^\infty {\hbox{A}}(s+u+x) {\hbox{A}}(s+v+x)\,
ds\cr
&=-{\hbox{A}}(u+x){\hbox{A}}(v+x);&(2.24)\cr}$$
hence 
$$ {{\partial }\over {\partial x}}\log\det \bigl( I-\kappa^2
\Gamma_{(x)}^2\bigr)=
-\kappa^2\bigl\langle (I-\kappa^2\Gamma_{(x)}^2)^{-1}{\hbox{A}}
(x+.),{\hbox{A}}(x+.)\bigr\rangle_{L^2(0,\infty )}.\eqno(2.25)$$
\indent By comparing the terms in the power series in $\kappa$, we deduce that

$$\kappa K(x,x)={{\partial }\over {\partial x}}\log\det
 (I-\kappa^2P_{(x,\infty )}WP_{(x,\infty )}).\eqno(2.26)$$
\rightline {$\square$}\par
\vskip.05in

\noindent {\bf Corollary 2.5.} {\sl Suppose further that} ${\hbox{A}}$ {\sl is an entire
function such that} $\int_0^\infty u\vert {\hbox{A}}(z+u)\vert^2\,du <\infty $ {\sl for 
each $z\in {\bf C}$, and let $\Gamma_{(z)}$ be the Hankel operator on
$L^2(0,\infty )$ that has kernel}
${\hbox{A}}(z+s+t)$. {\sl  Then
$${{d}\over {dz}} \log \det \bigl( I-\kappa^2 \Gamma^2_{(z)}\bigr)\eqno(2.27)$$
\noindent defines a meromorphic function on ${\bf C}$, extending $\kappa K(x,x)$.\/}\par
\vskip.05in
\noindent {\bf Proof.} By Morera's theorem,  $z\mapsto  \Gamma_{(z)}$ defines an entire 
function with values in $c^2$, and hence $ \det\bigl( I-\kappa^2
\Gamma^2_{(z)}\bigr)$ defines an entire function. The formula (2.27) defines
a holomorphic function,
except at those isolated points where the determinant vanishes, and these give rise to 
poles.\par
\rightline {$\square$}\par 
\vskip.05in

\indent In some cases $y={{d}\over {dx}}K(x,x)$ satisfies a Painlev\'e equation as in
[14, p 344]; that
is, $y''=F(y',y,x)$ where $F$ is rational in $y$ and $y'$, and analytic in $x$,
and such
that the only movable singularities of $y$ in ${\bf C}$ are poles. 
In particular, the bulk kernel gives rise to the $\sigma$ form of ${\hbox
{P}}_{\hbox{\smalletters V}}$, 
the hard edge ensemble gives rise to the ${\hbox{P}}_{\hbox{\smalletters III}}$ equation and the soft edge to 
${\hbox{P}}_{\hbox{\smalletters II}}$ as in [28, 29, 30, 12].\par
\indent The Painlev\'e ODE test asserts that every ordinary differential equation that 
arises from a partial differential equation via a Marchenko linear integral equation 
may be transformed to a Painlev\'e equation; see [1].\par
\vskip.05in
\indent A special feature of (2.10) is that $W$ is the
square of a self-adjoint Hankel operator $\Gamma$ on $L^2(0,\infty )$, and 
Proposition 2.2 gives an explicit construction of the symbol for $\Gamma$. Up to unitary
equivalence, this holds under some general spectral conditions which we list below. In section 3 we give 
a means for verifying (iii) for operators of the form (1.1).\par
\vskip.05in
\noindent {\bf Proposition 2.6.} {\sl Suppose that $W$ is a linear operator on $H$ such that\par
\indent (i) the nullspace of $W$ is either trivial or infinite-dimensional;\par
\indent (ii) $W$ is not invertible;\par
\indent (iii) $W$ is bounded and self-adjoint, and $W\geq 0$;\par
\indent (iv) $W$ has a simple discrete spectrum.\par
\noindent Then there exists a self-adjoint Hankel operator $\Gamma$ on $L^2(0,\infty )$
and a unitary
operator $U:H\rightarrow L^2(0,\infty )$ such that $W=U^*\Gamma^2 U$.}\par 
\vskip.05in

\noindent {\bf Proof.} Megretskii, Peller and Treil [21, p 257] have obtained sufficient conditions for a self-adjoint operator to be unitarily equivalent to 
the modulus of a Hankel operator. Under the more stringent condition (iv), 
their construction gives a Hilbert space $K$, a bounded linear operator
$X:K\rightarrow K$, and vectors $\xi ,\eta \in K$ 
such that the Hankel operator 
$$\Gamma f(t)=\int_0^\infty h(s+t)f(s)\, ds\qquad (f\in L^2(0,\infty ))\eqno(2.29)$$
\noindent with symbol $h(t) = \langle e^{tX}\xi ,\eta \rangle_K$ is unitarily
equivalent to $W^{1/2}$; thus $\Gamma$ is realized from a balanced
linear system in continuous time with one-dimensional input and output spaces.\par  
\vskip.05in
\noindent {\bf 3. Reproducing kernels and the bulk of the spectrum}\par
\indent In this section we recover the bulk kernel as a reproducing kernel, and show
more generally why operators of the form (1.1) are positive on weighted Hardy spaces. Let
$E$ be a meromorphic and zero-free function on ${\bf C}$ and let
$E^*(z)=\overline{E(\bar
z)}$, which has similar properties. We also introduce the meromorphic functions
$A(z)=(E(z)+E^*(z))/2$ and $B(z)=(E^*(z)-E(z))/(2i)$, which have $A(x)$ and $B(x)$
real for real $x$. \par
\indent Let $EH^2$ be the weighted Hardy space of meromorphic
functions $g$ on ${\bf C}_+$ such that $g/E$ belongs to the usual Hardy
space $H^2$, and with the
inner product 
$$\langle  g_1,g_2\rangle_{EH^2}=\langle g_1/E,
g_2/E\rangle_{H^2}=\int_{-\infty}^\infty 
{g_1(t)\bar g_2} (t){{dt}\over {\vert E(t)\vert^2}}.\eqno(3.1)$$
\noindent Similarly we can introduce $E^*H^2$. When $\zeta\in {\bf C}_+$ is not a 
pole of $E$, the linear functional $g\mapsto g(\zeta )$ is
bounded on $EH^2$, and hence given by $g(\zeta )=\langle g, k_\zeta\rangle_{EH^2}$,
where the reproducing kernel is  
$$k_\zeta (z)={{E(z)\overline{E(\zeta )}}\over {2\pi i(\bar\zeta -z)}}.\eqno(3.2)$$
\noindent We introduce $\Omega$ as the domain consisting
of points $z\in {\bf C}_+$, that are not poles of $E$ or $E^*$.
 Let $u(z)=E^*(z)/E(z)$, which is
meromorphic and unimodular on the real line, let 
$M_u:EH^2\rightarrow E^*H^2$ be the isometry $M_uf=uf$, and let $T_{\bar
u}:H^2\rightarrow H^2$ be the
Toeplitz operator $T_{\bar u}=R_+M_{\bar u}R_+$. \par
\vskip.05in
\noindent {\bf Theorem 3.1.} {\sl (i) The 
operator $W$ on $EH^2$ that has kernel
$$ W(z,w)={{E^*(z)\overline {E^*(w)}-E(z)\overline {E(w)}}\over {2\pi i(z-\bar
w)}}\qquad (z,w\in \Omega )\eqno(3.3)$$
\noindent compresses to an operator $EH^2\rightarrow EH^2$ that is unitarily equivalent 
to $\Gamma_{\bar u}^*\Gamma_{\bar u}$, where $\Gamma_{\bar u}:H^2\rightarrow 
\overline{H^2}$ is the Hankel operator $\Gamma_{\bar u} =R_-M_{u*}R_+$.}\par
\indent {\sl (ii) There exists a unique Hilbert space $H(W)$ of holomorphic functions on
$\Omega$ such that $W(z,w)$ is the reproducing kernel for $H(W)$.}\par   
\indent {\sl (iii) Suppose that $T_{\bar u}$ has a non-zero nullspace $K$. Then
$\Gamma_{\bar u}$
restricts to an isometry $K\rightarrow \overline{H^2}$.}\par
\vskip.05in
\noindent {\bf Proof. } (i) We write 
$$\eqalignno{\int_{-\infty}^\infty{{E^*(z)E(t)-E(z)\overline {E(t)}}\over 
{2\pi i(z-t)}}{{f(t)\, dt}\over {E(t)\overline{ E(t)}}}&=
{{E(z)}\over{2\pi i}}\int_{-\infty}^\infty {{f(t)/E(t)}\over {t-z}}\, dt\cr
 &\quad -{{E^*(z)}\over{2\pi i}}\int_{-\infty}^\infty
{{u^*(t)f(t)/E(t)}\over {t-z}}\, dt,&(3.4)\cr}$$
\noindent and hence by Cauchy's integral formula we have 
$$\eqalignno{Wf(z)&=E\bigl( f/E-M_uR_+M_{u^*} (f/E)\bigr)\cr
                 &=EM_uR_-M_{u^*} (f/E).&(3.5)\cr}$$
\noindent The map $V:EH^2\rightarrow H^2:$ $f\mapsto f/E$ is a unitary equivalence with
adjoint $V^*:g\mapsto Eg$, and $\Gamma_{\bar u}^*\Gamma_{\bar u} :H^2\rightarrow H^2$ reduces to 
$\Gamma_{\bar u}^*\Gamma_{\bar u}=R_+M_uR_-M_{u^*}R_+,$ so $\langle Wf,g\rangle_{EH^2}=\langle
V^*\Gamma_{\bar u}^*\Gamma_{\bar u} Vf,g\rangle_{H^2}$ for all $f,g\in EH^2$.\par
\indent (ii) By (i), $W$ is a positive operator on $EH^2$, so the
kernel $W(z,w)$ is of
positive type on $\Omega$; further, $z\mapsto W(z,w)$ and $w\mapsto W(z,\bar w)$ are
holomorphic on $\Omega$. Hence we can apply [2, Theorem 2.3.5] to obtain the Hilbert
space of holomorphic functions such that $W(z,w)$ is the reproducing kernel.\par 
\indent (iii) By [23, p 89], we have $T^*_{\bar
u}T_{\bar u}=I-\Gamma_{\bar u}^*\Gamma_{\bar u} ,$ which leads directly to the
identity $K=\{ f\in H^2: \Vert \Gamma_{\bar u}f\Vert =\Vert f\Vert\}.$\par
\rightline {$\square$}\par
\vskip.05in
\noindent {\bf Corollary 3.2.} {\sl Suppose that $u$ belongs to $H^\infty$
so that $E^*H^2$ is a closed linear subspace
of $EH^2$, and let $K=EH^2\ominus E^*H^2$ be the orthogonal complement of the range of
$M_u:EH^2\rightarrow EH^2$. Then $K$ equals $H(W)$ and has reproducing kernel}
$$K_w(z)={{A(z)B(\bar w)-B(z)A(\bar w)}\over{\pi (\bar w-z)}}\qquad(z,w\in
\Omega ).\eqno(3.6)$$
\vskip.05in
\noindent {\bf Proof.} First, one can check by calculation that
$$K_w(z)=W(z,w)={{E^*(z)\overline {E^*(w)}-E(z)\overline {E(w)}}\over {2\pi i(z-\bar w)}}.\eqno(3.7)$$
\noindent Then we observe that 
$$ {{E^*(z)\overline{E^*(w)}}\over {2\pi i(z-\bar w)}}
 =u(z){{E(z)\overline{E^*(w)}}\over {2\pi i(z-\bar w)}}\eqno(3.8)$$
\noindent lies in the range of $M_u$; so for $g\in K$ the proof of Theorem 3.1(i) 
simplifies to give
$$\bigl\langle g, K_w\bigr\rangle_{EH^2}=\bigl\langle g,
k_w\bigr\rangle_{EH^2}=g(w)\qquad (w\in \Omega ).$$
\rightline {$\square$}\par
\vskip.05in
\indent {\sl Bulk of the spectrum.} Thus when $u$ is an inner function we can identify $H(W)$ explicitly as the orthogonal complement of a
shift-invariant subspace of $EH^2$. In particular, by taking the entire function $E(z)=e^{-iaz}$, we find $u(z)=e^{2iaz}$ and 
the reproducing kernel 
for $K=EH^2\ominus E^*H^2$ to be  
$$K_w(z)={{\sin a(z-\bar w)}\over
{\pi (z-\bar w)}},\eqno(3.9)$$
\noindent as in the bulk kernel $B_{a/\pi}(z,w)$ of (1.3). Here we have 
$EH^2={\cal F}^*L^2[-a, \infty )$, and $\Psi_a={\cal F}^*\vert L^2[-a,a]$ gives a
unitary isomorphism $L^2[-a,a]\rightarrow K$ with $\Psi_a\Psi_a^*=B_{a/\pi}$. The Hankel operator $\Gamma_{\bar u}$ is isometric on $H^2\ominus
e^{2iax}H^2\simeq K.$\par   
\indent The Paley--Wiener theorem [17, p. 179]
characterizes $K$ as the space of functions $f\in L^2({\bf R})$ that are entire 
and of exponential type with 
$$\lim\sup_{y\rightarrow\pm \infty }\vert y\vert^{-1}\log \vert f(iy)\vert \leq 
a.\eqno(3.10)$$
\noindent Alternatively, we can characterize the subspaces by their scaling
properties. Let $(\delta_t)$ $(t\in {\bf R})$ be the unitary dilatation group on $L^2({\bf R})$ with 
$\delta_tf(x)=e^{t/2}f(e^{t} x)$. In [16], Katavolos and Power characterize the lattice of closed linear subspaces of $L^2$
that are simply invariant for both $S_s$ $(s\geq 0)$ and $\delta_s$ $(s\geq 0).$ \par 
\vskip.05in
\noindent {\bf Proposition 3.3.} {\sl The closed linear subspace $B_tL^2$ is simply invariant for $(\delta_s)$
$(s\leq 0)$, doubly invariant for $(\tau_s)$ $(s\in {\bf R})$ and
invariant under $J$. Conversely, if $\hat K$ is
any closed linear subspace of $L^2$ that is simply invariant for $(\delta_t)$
$(t\leq 0)$, doubly invariant for $(\tau_s)$ $(s\in {\bf R})$ and invariant under $J$, 
then $\hat K=B_aL^2$ for 
some $a>0$.}\par
\vskip.05in
\noindent {\bf Proof.} We have $\delta_{-s}={\cal F}^*\delta_s{\cal F}$ and $\tau_s={\cal
F}^* S_{-s}{\cal F}$, so we shall characterize the subspaces $L^2[-\pi t, \pi t]$
under the operation of $\delta_s$, $S_s$ and $J$. Now $L^2[-t\pi ,\pi t]$ is clearly doubly invariant for
$(S_s)$ $(s\in {\bf R})$, and $\delta_s L^2[-\pi t, \pi t]=L^2[-\pi te^{-s}, \pi
te^{-s}]$; so $L^2[-t\pi ,\pi t]$ is simply invariant for $(\delta_s )$ $(s\geq 0)$.
Conversely, all closed linear subspaces $\hat K$ of $L^2$ that are simply invariant under
$(\delta_s )$ $(s\geq 0)$ and doubly invariant under $(S_s)$ $(s\in {\bf R})$ have the
form $\hat K=L^2(-a,b)$ for some $a, b\in {\bf R}\cup\{\infty \}$ by a simple case of
Beurling's theorem. When $\hat K$ is additionally invariant under
$J$, we need to have $a=b$; hence $\hat K=L^2[-a,a]$.\par    

\vskip.1in

\noindent {\bf 4. Hard-edge operators and Sonine spaces}\par
\indent In this section we consider the operators for the hard edge case and
associated subspaces. Let
$J_\nu$ be the Bessel function of the first kind for real $\nu >-1/2$, and let  
$$h(z)=\sum_{k=0}^\infty {{(-1)^k(1+2ik)z^k}\over {2^{\nu +2k}\Gamma (\nu
+k+1)k!}}=z^{-\nu /2}J_\nu(\sqrt z)+2iz{{d}\over {dz}}\Bigl( z^{-\nu /2}J_\nu(\sqrt
z)\Bigr)\eqno(4.1)$$
\noindent which is entire and of order $1/2$ as in [11, p 190]. Then $E(z)=1/h(z)$ is a meromorphic
function, with no zeros, such that 
$${{E^*(z)\overline{E^*(w)}-E(z)\overline {E(w)}}\over {2\pi i(z-\bar w)}}$$
$$=\Bigl({{J_\nu (z^{1/2})w^{1/2}J_\nu'(\bar w^{1/2})-z^{1/2}J_\nu'(z^{1/2})J_\nu (\bar
w^{1/2})}\over {\pi (z-\bar
w)}}\Bigr)\Bigl({{E(z)E^*(z)\overline{E^*(w)}\overline{E(w)}}\over {z^{\nu
/2}\bar w^{\nu/2}}}\Bigr).\eqno(4.2)$$
\noindent We recognise the first factor on the right-hand side from (1.8), and the left-hand side from
(4.4); but
Corollary 3.2 does not apply directly to $E^*(z)/E(z)$; so we introduce
operators that correspond to these kernels indirectly by means of
the Hankel transform as in [25, p 298]. The Hankel transform 
of $f\in L^2(xdx; (0,\infty ))$ is
$${\cal H}_\nu f(x)=\int_0^\infty J_\nu(xy) f(y)\, ydy.\eqno(4.3)$$
\indent On $L^2(xdx; (0,\infty ))$ we introduce the unitary dilatation 
group $(\tilde \delta_t)$ by
$\tilde \delta_tg(x)=e^tg(e^tx)$ and the unitary operator $U:L^2(xdx, (0,\infty ))\rightarrow
L^2({\bf R})$ by $Ug(\xi )=e^{-\xi }g(e^{-\xi})$ such that $U^*\tau_t
U=\tilde\delta_t$.\par 
\vskip.05in
\noindent {\bf Lemma 4.1.} {\sl Let $G_\ell$ be the integral operator on $L^2({\bf R})$
that has kernel function
$$e^{-\ell -\xi-\eta} J_\nu (e^{-\ell -\xi-\eta}).\eqno(4.4)$$
\indent Then $G_\ell$ is a self-adjoint and unitary operator such that
$G_\ell^2=I$, and $G_\ell\tau_t=\tau_{-t}G_\ell$.\/}\par
\vskip.05in
\noindent {\bf Proof.} From the shape of the integral kernel, the identity $G_\ell
U=\tau_{-\ell} U{\cal H}_\nu$ is evident. Further, Hankel's inversion formula leads to the
identity ${\cal H}_\nu^2=I$, whence to 
$$G_\ell UU^*G_\ell =\tau_{-\ell} U{\cal H}_\nu {\cal H}_\nu U^* \tau_\ell
=I.\eqno(4.5)$$
\indent The identity (4.4) is evident from the definitions, and by (4.5) is
equivalent to the scaling property 
${\cal H}_\nu \tilde\delta_t=\tilde\delta_{-t}{\cal H}_\nu$ of the Hankel transform
as in [25, p. 299].\par
\rightline {$\square$}\par
\vskip.05in
\indent The following result on position of subspaces contrasts with Corollary 3.2. Here $\Gamma$
denotes Euler's gamma function.\par
\vskip.05in
\noindent {\bf Theorem 4.2.} {\sl (i) The operator $Q_\ell =G_\ell P_+G_\ell$ on
$L^2({\bf R})$ is an
orthogonal projection.\par
\indent (ii) The range of
${\cal F}Q_\ell {\cal F}^*$ equals $e^{i\ell x}u_\nu H^2,$ where the
meromorphic function}\par
$$u_\nu (z)=2^{iz}{{\Gamma ((1+\nu +iz)/2)}\over{\Gamma ((1+\nu -iz)/2)}}\eqno(4.6)$$
\noindent {\sl is holomorphic on $\{z:\Im z<0\}$, and unimodular and continuous on ${\bf R}$.\par
\indent (iii) Whereas $u_\nu^*H^2\cap H^2=\{ 0\},$ for $\nu >0$ the subspace $K=(u_\nu \overline{H^2})\cap H^2$ is
non-zero, and $\Gamma_{\bar u_\nu}:H^2\rightarrow
\overline {H^2}$ restricts to an isometry $K\rightarrow
\overline {H^2}$.}\par 

\vskip.05in
\noindent {\bf Proof.} (i) This follows directly from the Lemma.\par
\indent (ii) Our aim is to show that the range of the
orthogonal projection ${\cal F}^*Q_\ell {\cal F}$ is simply invariant under the
$S_\lambda $ for $\lambda >0$. By Plancherel's theorem we have
$$S_\lambda {\cal F}Q_\ell L^2={\cal F}\tau_{-\lambda}G_\ell P_+L^2={\cal
F}G_0\tau_{\lambda +\ell} P_+L^2,\eqno(4.7)$$
\noindent where $\tau_{\lambda +\ell}P_+L^2=L^2(\lambda +\ell , \infty )\subseteq
L^2(\ell , \infty )$ and $\cap_{\lambda >0}L^2(\lambda , \infty )=0$. Consequently by
Beurling's theorem, there exists a unimodular and measurable function $u_\nu$ such that 
${\cal F} Q_0 u_\nu L^2=u_\nu H^2$, and $u_\nu$ is unique up to a unimodular
constant factor. One can easily deduce that ${\cal F} Q_\ell
L^2=e^{i\ell x}u_\nu H^2$.\par
\indent The Fourier conjugate of $Q_\ell$ is ${\cal F}Q_\ell {\cal F}^*=
{\cal F}G_\ell {\cal F}^*{\cal F}P_+{\cal F}^*{\cal F}G_\ell {\cal F}^*,$ wherein we recognise ${\cal F}P_+{\cal F}^*$ as 
$R_-:L^2\rightarrow \overline{H^2}$. To determine the range of ${\cal F}Q_\ell {\cal
F}^*$, or
equivalently the subspace ${\cal F}G_\ell L^2(0,\infty )$, we write 
$${\cal F}G_\ell f(x)=\int_{-\infty}^\infty e^{-ix\xi} \int_0^\infty e^{-\ell -\xi
-\eta} J_\nu (e^{-\ell -\xi-\eta}) f(\eta )d\eta {{d\xi}\over {\sqrt{2\pi}}}$$
\noindent for $f\in L^2(0,\infty )$, and then reduce this integral by simple
transformations to
$${\cal F}G_\ell f(x)=e^{ix\ell} {\cal F^*} f(x)\int_{-\infty}^\infty e^{-(1+\nu
+ix)\xi}
e^{\nu \xi }J_\nu (e^{-\xi })\, d\xi .\eqno(4.8)$$
\noindent The substitution $y=e^{-\xi}$ reduces the final integral in (4.8) to a standard Mellin
transform [25, p. 263], and we identify $u_\nu$ from 
$${\cal F}G_\ell f(x)=e^{ix\ell} {{2^{ix}\Gamma ((1+\nu +ix )/2)}\over {\Gamma ((1+\nu
-ix )/2)}}{\cal F}^* f(x).\eqno(4.9)$$
\indent (iii) Let $E_\nu (z)=e^{-iz\log\sqrt 2}\Gamma ((1+\nu -iz)/2)$ so that 
$E_\nu $ is meromorphic and zero-free with simple poles at $-i-\nu i-2ki$
for $k=0,1,\dots $, and $u_\nu(z)=E_\nu^*(z)/E_\nu(z)$ has simple zeros at 
$z_k=-i-\nu i-2ki$ for $k=0, 1, \dots $  and simple poles at 
$i+\nu i+2ki$ for $k=0, 1, \dots $. The function $u_\nu (z)$ is holomorphic in the
lower half plane, but does not define a bounded holomorphic function on $\{ z:\Im
z<0\}$ since the series  $\sum_{k=0}^\infty \Im z_k /(1+\vert z_k\vert^2)$ diverges,
violating Blaschke's condition for the zeros of a non-trivial function in $H^\infty$
or $H^2$ as in [17, p. 92]. Hence the equations $h_1(z)=u_\nu^*(z)h_2(z)$ with
$h_1,h_2\in H^2$ has only the trivial solution $h_1=h_2=0$; so $u^*_\nu H^2\cap H^2=0.$ 
Note that $\log \vert u^*_\nu (z)\vert$ is subharmonic on the
${\bf C}_+$, but is not the Poisson integral of a measure on ${\bf R}$.\par
\indent We take $a>0$ and $\nu +1/2>\lambda >1/2 $, and let
$$f(x)=a^{\nu -\lambda +3/2}x^{1/2-\nu }(x^2-a^2)^{(\lambda -1)/2}J_{\lambda -1}
\bigl(a\sqrt{x^2-a^2}\bigr){\bf I}_{(a, \infty )}(x),$$
\noindent with Hankel transform
$$g(t)=t^{1/2}{\cal H}_\nu \bigl( x^{-1/2} f(x); t\bigr).$$
\noindent Then by a result of Sonine [6 p. 301 , 24 p. 75, 26 p. 38], both $f$ and $g$ are supported on
$(a,\infty )$ and  we have 
$$\int_a^\infty g(t) t^{-1/2+ix}\, dt=u_\nu (x) \int_a^\infty f(t)t^{-1/2-ix}\, dt\qquad
(x\in {\bf R}).
\eqno(4.10)$$
\noindent Hence when $a=1$ there exist non-zero functions $h_1,h_2\in H^2$ such that
$h_2(x)=u_\nu (x)h^*_1(x)$, so $h_2\in
u_\nu\overline{H^2}.$ Now we apply Theorem 3.1(iii) to deduce that  
$\Gamma_{\bar u_\nu}\vert H^2\cap u_\nu\overline{H^2}$ is an isometry.\par
\rightline{$\square$}\par 
\vskip.05in
\noindent {\bf Proposition 4.3.} {\sl (i) The Hankel operator $\Phi_\ell =P_+G_\ell P_+$
on $L^2(0,\infty )$ has $\Phi_\ell^2 =P_+Q_\ell P_+$.\par
\indent (ii) The operator $\Phi_\ell$ on $L^2(0,\infty )$ is Hilbert--Schmidt, and 
each non-zero\par
\noindent $f\in L^2(xdx, (0,1))$ such that
$$\lambda f(x)=\int_0^1 J_\nu (\sqrt{sxy})f(y)\, dy\eqno(4.11)$$
\noindent corresponds to an eigenfunction $g\in
L^2(0,\infty)$ of $\Phi_\ell$ with eigenvalue ${{1}\over {2}}\lambda \sqrt s$.\par
\indent (iii) The kernel of $Q_\ell$ as an integral operator on $L^2({\bf R})$ is} 
$${{e^{-\ell -\xi}J_\nu (e^{-\ell -\xi})e^{-2\ell-2\eta}J_\nu'(e^{-\ell -\eta })-
e^{-2\ell-2\xi}J_\nu'(e^{-\ell -\xi})e^{-\ell -\eta}J_\nu (e^{-\ell -\eta})}\over
{e^{-2\ell-2\xi}-e^{-2\ell-2\eta}}}.\eqno(4.12)$$
\indent {\sl (iv) $\det (I-zF^{0,a})=\det (I-z\Phi_{(\alpha )}^2)$ for
$\alpha =-(1/2)\log a$ and $a>0$. }\par
\vskip.05in
\noindent {\bf Proof.} (i) For $t>0$ we have the Hankel condition
$\Phi_\ell\tau_t=\tau_{t}^*\Phi_\ell,$ where here $(\tau_t)_{t>0}$ denotes the
semigroup of translation operators on $L^2(0,\infty )$. Then one uses Theorem 4.2(i).\par
\indent (ii) The kernel function is clearly symmetric, real-valued and square
integrable, since
$$\int_0^\infty\!\!\!\int_0^\infty e^{-2(\ell +\eta +\xi)}J_\nu (e^{-(\ell +\eta
+\xi)})^2\, d\xi d\eta =\int_0^\infty ue^{-2\ell -2u}J_\nu (e^{-\ell -u})^2\,
du<\infty\eqno(4.13)$$
\noindent due to the asymptotic formula $J_\nu (x)\asymp x^\nu/ \Gamma (\nu +1)$ as
$x\rightarrow 0+.$ Hence $\Phi_\ell$ gives a self-adjoint operator of Hilbert--Schmidt type. The
operator $U$ restricts to a unitary $L^2(xdx; (0,1))\rightarrow L^2(0, \infty )$, and
under this transformation the
eigenfunction equations correspond via $g(\xi )=e^{-\xi} f(e^{-2\xi }).$\par
\indent (iii) We use the method of proof of Proposition 2.2 to verify the stated
formula for $Q_\ell =G_\ell P_+G_\ell$, which is 
essentially the square of a
self-adjoint Hankel operator. 
With $A(\xi )=e^{-\xi }J_\nu (e^{-\xi})$ and $B(\xi
)=e^{-2\xi}J_\nu'(e^{-\xi})$, we have 
$${{d}\over{d\xi}}\left[\matrix{A\cr B\cr}\right] = 
             \left[\matrix{-1&-1\cr (e^{-2\xi} -\nu^2)&-1\cr}\right] 
\left[\matrix{A\cr B\cr}\right]\eqno(4.14)$$
\noindent where 
$$\left[\matrix{0&-1\cr 1 &0\cr}\right]\left[\matrix{-1&-1\cr (e^{-2\xi} 
-\nu^2)&-1\cr}\right]+\left[\matrix{-1& (e^{-2\eta} 
-\nu^2)\cr -1&-1\cr}\right]\left[\matrix{0&-1\cr 1 &0\cr}\right]$$
$$=\left[\matrix{e^{-2\eta}-e^{-2\xi}&0\cr 0 &0\cr}\right]+
\left[\matrix{0&-2\cr 2 &0\cr}\right],\eqno(4.15)$$
\noindent hence 
$$\Bigl( {{\partial}\over {\partial \xi}} + {{\partial}\over {\partial \eta}}\Bigr) 
{{A(\xi )B(\eta )-A(\eta )B(\xi )}\over {e^{-2\xi}-e^{-2\eta}}} = 
\Bigl( {{\partial}\over {\partial \xi}} + {{\partial}\over {\partial \eta}}\Bigr)
\int_0^\infty A(\xi +u)A(\eta +u)\, du.$$
\noindent Hence we can obtain the stated identity by following the proof of Proposition
 2.2. Alternatively, one can transform a formula in [25, p. 303]. \par
\indent (iv) The unitary equivalence between $L^2((0,1), dx)$ and
$L^2(0,\infty )$ involves $g(x)\mapsto \sqrt 2e^{-\xi}g(e^{-2\xi}),$ so 
$F^{(0,1)}$ is unitarily equivalent to the operator that has kernel
$$2e^{-\xi-\eta }F^{(0,1)}(e^{-2\xi }, e^{-2\eta})= {{e^{-\xi}J_\nu
(e^{-\xi})e^{-2\eta}J_\nu'(e^{-\eta})-e^{-2\xi}J_\nu'(e^{-\xi})
e^{-\eta}J_\nu(e^{-\eta})}\over{e^{-2\xi}-e^{-2\eta}}},$$
\noindent which we recognise as the kernel of
$\Phi_{(0)}^2$. Comparing the spectra of the compressions to $L^2(0,a)$
and $L^2(\alpha ,\infty )$, we deduce that 
$$\det (I-zF^{0,a})=\det(I-zP_{(\alpha ,\infty
)}\Phi_{(0)}^2P_{(\alpha ,\infty )})=\det (I-z\Phi_{(0)}P_{(\alpha ,\infty
)}\Phi_{(0)}).\eqno(4.16)$$
\noindent Finally, $\Phi_{(0)}P_{(\alpha
,\infty )}\Phi_{(0)}$ equals $\Phi_{(\alpha )}^2$ since they both have
kernel
$$\int_\alpha^\infty e^{-\xi-u}J_\nu (e^{-\xi -u})e^{-\eta -u}J_\nu
(e^{-\eta -u})\, du.\eqno(4.17)$$  
\rightline{$\square$}\par    
\vskip.05in

\noindent {\bf Theorem 4.4.} {\sl Let $T$ be the operator 
$$Tf(\xi )=-{{\partial}\over {\partial \xi}}\Bigl( e^{2\xi}{{\partial}f\over 
{\partial \xi}}\Bigr) +(\nu^2-1)e^{2\xi} f(\xi ).\eqno(4.18)$$
\indent (i) Then $T$ is an essentially self-adjoint and positive operator on
$C_c^\infty({\bf R})$ in $L^2({\bf
R})$, so that $V_t=e^{-it\ell}T^{-it/2}$ $(t\in {\bf R})$ defines a $C_0$ group of unitary
operators on $L^2({\bf R})$.\par
\indent (ii) The unitary groups $(V_s)_{s\in {\bf R}}$ and $(\tau_t )_{t\in {\bf R}}$
satisfy $V_s\tau_t =e^{ist} \tau_tV_s$ for 
$s,t\in {\bf R}.$ \par
\indent (iii) The subspace $Q_\ell
L^2$ is doubly invariant for $(V_s)$ with $s\in {\bf R}$ and simply invariant
for $(\tau_{-t})$ for $t\geq 0$. Conversely, if $K$ is a non-trivial closed
linear subspace of $L^2$ that is simply invariant for $\tau_{-t}$ $(t\geq 0)$ and
doubly invariant for $V_s$ $(s\in {\bf R})$, then $K=Q_\alpha L^2$ for some real $\alpha$.}\par

\vskip.05in
\noindent {\bf Proof.} (i) The simplest way of proving that the operator $T$ is self-adjoint
is to compute its spectral resolution. By simple transformations of the Bessel
equation [14, p 171], we have  
$$-e^{2\xi}\Bigl( {{\partial^2}\over {\partial \xi^2}}+2{{\partial}\over {\partial \xi}}
+\nu^2-1\Bigr) \bigl(e^{-\xi -\ell-\eta }J_\nu (e^{-\xi -\ell-\eta })\bigr)= e^{-2\ell-2\eta
}\bigl(e^{-\xi -\ell-\eta }J_\nu (e^{-\xi -\ell-\eta })\bigr),\eqno(4.19)$$
\noindent so that $e^{-\xi -\ell-\eta }J_\nu (e^{-\xi -\ell-\eta })$ is an eigenfunction
of $T$ corresponding to the eigenvalue $e^{-2\ell-2\eta}>0$. By Hankel's inversion
theorem [25, p 299], the functions $\lambda yJ_\nu (\lambda xy)$ give a complete spectral family in 
$L^2(xdx;
(0,\infty ))$, and the unitary transformation $U$ takes $\lambda yJ_\nu (\lambda xy)$ to 
$e^{-\xi-\ell -\eta }J_\nu (e^{-\xi -\ell-\eta})$ after an obvious change of variable.
By Stone's theorem, $(-i/2)\log T$ generates a $C_0$ unitary group
$T^{-is/2}$.\par
\indent (ii) We have
$$\eqalignno{ V_sG_\ell f(\xi )&=e^{-is\ell}T^{-is/2}\int_{-\infty}^\infty e^{-\xi -\eta
-\ell} J_\nu ( e^{-\xi -\ell-\eta})f(\eta)\, d\eta\cr
 &=\int_{-\infty}^\infty e^{is\eta}e^{-\xi -\ell-\eta }J_\nu (e^{-\xi
-\ell-\eta})f(\eta )\, d\eta\cr
&=G_\ell S_sf(\xi );&(4.20)\cr}$$
\noindent hence $G_\ell V_sG_\ell =S_s.$ When we conjugate the familiar Weyl--von
Neumann relation $\tau_{-t}S_s=e^{ist}S_s\tau_{-t}$ by $G_\ell\otimes G_\ell^*$ we 
obtain $ G_\ell\tau_{-t}G_\ell G_\ell S_sG_\ell=
e^{ist}G_\ell S_sG_\ell G_\ell\tau_{-t}G_\ell$ or $\tau_t V_s=e^{ist} V_s \tau_t .$\par
\indent (iii) From earlier relations, we have 
$$V_sQ_\ell =V_sG_\ell P_+G_\ell =G_\ell S_sP_+=G_\ell P_+S_s=G_\ell P_+G_\ell G_\ell
S_s=Q_\ell G_\ell S_s,\eqno(4.21)$$
\noindent which shows that the range of $Q_\ell$ is mapped onto itself by $V_s$; further
$$\tau_{-t}Q_\ell =\tau_{-t}G_\ell P_+G_\ell =G_\ell \tau_tP_+G_\ell =
G_\ell P_{[t,\infty )}\tau_tG_\ell =G_\ell P_{[t,\infty )}G_\ell\tau_{-t},\eqno(4.22)$$
has range contained in the range of $Q_\ell$ for $t>0,$ so $Q_\ell L^2$ is simply
invariant.\par 
\indent To obtain the converse, we consider the Fourier transforms of the groups.
On\par
\noindent $e^{i\ell x}u_\nu H^2$, the unitary semigroups operate as 
$$\hat V_s={\cal F}V_s{\cal F}^*:f (x)\mapsto e^{i\ell s} u_\nu (x)u_\nu (s-x)f(x-s)\qquad
(s\in {\bf R});\eqno(4.23)$$
\noindent ${\cal F}\tau_{-t}{\cal F}^* =S_t$. To
verify (4.23), we recall the flip map $J$ by $Jf(x)=f(-x)$, and observe that ${\cal F}
{\cal F}=J$, and ${\cal F}^*{\cal F}^*=J$. We have
$${\cal F}V_s{\cal F}^*={\cal F}G_\ell S_sG_\ell {\cal F}^*=S_\ell M_{u_\nu} {\cal F}^*S_s{\cal
F}{\cal F}^*G_\ell {\cal F}^*\eqno(4.24)$$
\noindent so that ${\cal F}V_s{\cal F}^*=S_\ell M_{u_\nu}\tau_s JS_\ell M_{u_\nu} J.$ Using the 
Weyl--von Neumann relation for $\tau_s$ and $S_\ell$, one can
easily simplify this expression to obtain $\hat V_s={\cal F}V_s{\cal F}^*=
e^{i\ell s}M_{u_\nu}N_s\tau_s,$ where $N_sf(x)=u_\nu (s-x)f(x).$ The functions $u_\nu$ satisfy $u_\nu
(-x)u_\nu (x)=1$ and 
$$e^{i\ell s} u_\nu (x)u_\nu (s-x)=e^{i\ell s}2^{is} {{\Gamma ((1+\nu +ix)/2)
\Gamma ((1+\nu +is-ix)/2)}\over {\Gamma ((1+\nu -ix)/2)\Gamma ((1+\nu
+ix-is)/2)}}.\eqno(4.25)$$   

\indent Suppose that $K$ is such an invariant subspace. Then by Beurling's theorem,
there exists a unimodular and measurable function $w$ such that 
${\cal F}K=wH^2$; further, this $w$ is uniquely determined up to a unimodular constant multiple.
We apply ${\cal F}V_s{\cal F}^*$ to this identity, and deduce by double invariance and (4.23)
that 
$$\{ e^{i\ell s} u_\nu (x)u_\nu (s-x)w(x-s)f(x-s): f\in H^2\} =wH^2;\eqno(4.26)$$
\noindent so that, 
$$ e^{i\ell s} u_\nu (x)u_\nu (s-x)w(x-s)=c(s)w(x)\qquad (s\in {\bf R})$$
\noindent holds for some $c(s)$. We re-arrange this to $u_\nu (s-x)w(x-s)=u_\nu
(-x)w(x) e^{-i\ell s} c(s)$, then solve to obtain $w(x)=e^{i\alpha x}u_\nu (x)$ for
some $\alpha\in {\bf R}$. Hence ${\cal F}K=e^{i\alpha x}u_\nu
(x)H^2$, so 
$K=Q_\alpha L^2$ by Theorem 4.2(ii).\par
\rightline{$\square$}\par
\vskip.05in 
\noindent {\bf Corollary 4.5.} {\sl The unitary groups $(\hat V_s)$ $(s\in {\bf R})$ and
$(S_t)$ $(t\in {\bf R})$ form a Weyl pair. The closed linear subspaces that are simply
invariant for $\hat V_s$ $(s\geq 0)$ and $S_t$ $(t\geq 0)$ have the form}
$$K_{\alpha , \beta}=\{ e^{i\alpha x-i\beta x^2} u_\nu (x)f(x): f\in H^2\}\qquad
(\alpha\in {\bf R}, \beta >0).\eqno(4.27)$$
\vskip.05in
\noindent {\bf Proof.} We shall check simple invariance of $K_{\alpha , \beta}$ under 
$\hat V_s$ $(s\in {\bf R})$, since the other statements follow easily from the proof of
Theorem 4.4 and the Katavolos--Power theorem. We have 
$$\hat V_s:e^{i\alpha x-i\beta x^2} u_\nu (x)f(x)\mapsto e^{i\ell s-i\alpha
s-i\beta s^2} e^{i\alpha x-i\beta x^2} u_\nu (x)e^{2i\beta sx}f(x-s),\eqno(4.28)$$
\noindent where $e^{2i\beta sx}f(x-s)\in S_{2\beta s}H^2\subset H^2$ for $s,\beta
>0$.\par
\vskip.1in

\noindent {\bf 5 Soft edge operators and the Airy group}\par
\noindent With $D=-i{{\partial }\over {\partial
x}}$ 
the Airy group $e^{itD^3}$ is a $C_0$ group of unitary operators, as defined by
$$e^{itD^3}f(x)={{1}\over{\sqrt{2\pi}}}\int_{-\infty}^\infty e^{it\xi^3 +i\xi
x}\, {\cal F}f(\xi )\, d\xi .\eqno(5.1)$$ 
In this section we shall consider how the Airy group is related to the kernel function of
(1.11). Here $J_t$ denotes the operator $e^{itD^3}J$ on $L^2({\bf R})$, not a Bessel function, and 
 we shall use a subscript $t$ to indicate scaling of
the space variables $x$ and $y$ with respect to time $t$.\par
\vskip.05in
\noindent {\bf Lemma 5.1.} {\sl The operator $J_t=e^{itD^3}J$ is 
self-adjoint with $J_t^2=I$, and $J_t$ as an integral operator on $L^2({\bf R})$
has kernel}
$${{1}\over {(3t)^{1/3}}}{\hbox{Ai}}\Bigl( {{x+y}\over {(3t)^{1/3}}}\Bigr)
.\eqno(5.2)$$
\vskip.05in
\noindent {\bf Proof.} For a compactly supported and smooth function $f$ we have
$$Je^{itD^3} Jf(x)={{1}\over {\sqrt{2\pi}}}\int_{-\infty}^\infty 
e^{it\xi^3 -i\xi x} {\cal F}f(\xi )\, d\xi=e^{-itD^3} f(x),\eqno(5.3)$$
\noindent so $J_t^2=I$. Further, the kernel of $J_t$ is given by 
$$\eqalignno{ e^{itD^3}Jf(x)&={{1}\over{2\pi}}\int_{-\infty}^\infty 
e^{it\xi^3+i\xi x}\int_{-\infty}^\infty e^{i\xi y}f(y)\, dy \, d\xi\cr
&=\int_{-\infty}^\infty \Bigl\{ {{1}\over {2\pi}}\int_{-\infty}^\infty
e^{i\xi^3t +i\xi (x+y)}\, d\xi\Bigr\} f(y)\, dy\cr
&=\int_{-\infty}^\infty {{1}\over {(3t)^{1/3}}}{\hbox{Ai}}\Bigl( {{x+y}\over 
{(3t)^{1/3}} }\Bigr) f(y)\, dy.&(5.4)\cr}$$ 
\noindent Since the Airy function is real valued, it also follows that
$J_t$ is self-adjoint.\par
\rightline {$\square$}\par
\indent Most of the next result is essentially contained in [28, Lemma
2], but we
 include a proof for completeness.\par
\vskip.05in
\noindent {\bf Proposition 5.2.} {\sl (i) The operator 
$$W_t=e^{itD^3}P_-e^{-itD^3}=J_tP_+J_t\eqno(5.5)$$
\noindent on $L^2({\bf R})$ is an
orthogonal projection and the range of ${\cal F}W_t{\cal F}^*$ equals 
$e^{it\xi^3}H^2.$\par
\indent (ii) The Hankel operator $\Gamma_{0,t}=P_+J_tP_+$ has square
$\Gamma_{0,t}^2=P_+W_tP_+$.\par
\indent (iii) The kernel of $W_t$ as an integral operator on $L^2({\bf R})$ is
$$W_t(x,y)= {{1}\over {(3t)^{1/3}}}W_{1/3}\Bigl({{x}\over {(3t)^{1/3}}},
{{y}\over {(3t)^{1/3}}}\Bigr)\eqno(5.6)$$
\noindent where $W_{1/3}$ is the Airy kernel as in (1.11).\/}\par
\vskip.05in
\noindent {\bf Proof.} (i) By Lemma 5.1, we have
 $W_t^2=J_tP_+J_t^2P_+J_t=J_tP_+J_t=W_t$, so that $W_t$ is a
projection; further $W_t^* =W_t$. The range of $W_t$ equals $\{ W_tk:k\in L^2({\bf R})\}$, or
equivalently the range of $J_tP_+$.\par
\indent If $f\in L^2({\bf R}_+)$, then
${\cal F}f(\xi )=\bar G(\xi )$, where $G\in H^2.$ Since $e^{-itD^3}$ is unitary, we have
$W_tL^2=e^{itD^3}P_-e^{-itD^3}L^2=e^{itD^3}P_-L^2$, and hence the image
of $WL^2$ under the Fourier transform ${\cal F}$ is ${\cal F}W_tL^2=\{
e^{it\xi^3}F(\xi ):F\in H^2\}$.\par
\indent (ii) We have $\Gamma_{0,t}=P_+e^{itD^3}JP_+$ and hence
$$\eqalignno{\Gamma_{0,t}^2&=P_+e^{itD^3}JP_+e^{itD^3}JP_+\cr
&=P_+e^{itD^3} JP_+JJe^{itD^3}JP_+\cr
&=P_+e^{itD^3}P_-e^{-itD^3}P_+=P_+W_tP_+.&(5.7)\cr}$$
\indent (iii) It also follows from the Lemma that
the kernel function is
$$W_t(x,y)={{1}\over {(3t)^{2/3}}}\int_0^\infty {\hbox{Ai}}\Bigl( {{x+u}\over 
{(3t)^{1/3}}}\Bigr) {\hbox{Ai}}\Bigl({{u+y}\over{(3t)^{1/3}}}\Bigr)
du,\eqno(5.8)$$
\noindent a formula which reduces to (5.6) and (5.7) on account of the
identity 
$$\int_0^\infty
{\hbox{Ai}}(x+u){\hbox{Ai}}(u+y)\,
du={{{\hbox{Ai}}(x){\hbox{Ai}}'(y)-{\hbox{Ai}}'(x){\hbox{Ai}}(y)}
\over {x-y}}.\eqno(5.9)$$
\noindent This formula is presented by Tracy and Widom in [{28}], and 
follows from Proposition 2.2 since $A(x)={\hbox{Ai}}(x)$ satisfies
$${{d}\over{dx}}\left[\matrix{A\cr B\cr}\right]=\left[\matrix{0&1\cr x &0\cr}\right]
\left[\matrix{A\cr B\cr}\right].$$

\rightline {$\square$}\par

\indent Evidently $W_0=P_-$, and the relative positions of the ranges of $P_-$ and $W_t$ are
described in the following Proposition, with different conclusions from Theorem 4.2.\par
\noindent {\bf Definition.} [9] A function $G\in H^2$ is said to be
cyclic (for the backward shifts) when ${\hbox{span}}\{ S_t^*G; t>0\}$
is dense in $H^2$. Likewise, $f\in L^2({\bf R}_+)$ is cyclic when 
${\hbox{span}}\{\tau_t^*f: t>0\}$ is dense in $L^2({\bf R}_+)$;
 $g\in L^2({\bf R}_-)$ is cyclic when ${\hbox{span}}\{\tau_t^*g: t<0\}$ is dense in $L^2({\bf R}_-)$.\par
\vskip.05in

\noindent {\bf Proposition 5.3.} {\sl (i) For each $t\neq 0$, the subspaces
$W_tL^2 \cap L^2({\bf R}_-)$ and $(W_tL^2)^\perp \cap L^2({\bf R}_+)$ equal $\{0\}$; while 
any non-zero vector in 
$W_tL^2 \cap L^2({\bf R}_+)$ or $(W_tL^2)^\perp  \cap L^2({\bf R}_-)$ is cyclic.}\par
\indent {\sl (ii) 
For each $t>0$, the operator 
$W_t$ on $L^2({\bf R}_-)\oplus L^2({\bf
R}_+)$ has block matrix form}
$$\left[\matrix{P_-W_tP_-& P_-W_tP_+\cr P_+W_tP_-&P_+W_tP_+\cr}\right]
\in \left[\matrix {{\hbox{B}}&c^2
 \cr c^2& c^1\cr}\right].\eqno(5.10)$$
\indent {\sl (iii) For any real $t$, the operators $P_+W_tP_-$
and $P_-W_tP_+$ are Hilbert--Schmidt.\/}\par
\vskip.05in
\noindent {\bf Proof.} (i) First we check that $W_tL^2\cap L^2({\bf R}_-)=0$, or 
equivalently by Proposition 5.2(i) that $e^{it\xi^3}H^2\cap H^2=0.$ Suppose that 
$F,G\in H^2$ are non-zero and
satisfy $e^{it\xi^3}F(\xi )=G(\xi )$ for almost all $\xi \in {\bf R}$. Then $K(\zeta )=e^{it\zeta^3}F(\zeta )-G(\zeta )$ is a holomorphic
function with zero boundary values
at almost all points of ${\bf R}$; so by 
the Lusin--Privalov theorem, $K(\zeta )$ is identically zero on ${\bf C}_+$.  
 Now
by Szeg\"o's Theorem [{17}, page 108], the integrals
$$\int_{-\infty}^\infty {{\log \vert F(\xi +i\eta )\vert
}\over {1+\xi^2}} d\xi  \quad {\hbox{and}}\quad\int_{-\infty}^\infty 
{{\log \vert G(\xi +i\eta )\vert
}\over {1+\xi^2}} d\xi \eqno(5.11)$$    
\noindent converge. But this contradicts the identity 
$e^{it\zeta^3}F(\zeta )=G(\zeta )$, since 
$$t\int_{-\infty}^\infty {{\Im (\xi +i\eta )^3}\over {1+\xi^2}}
d\xi\eqno(5.12)$$
\noindent diverges for $\eta ,t>0$; so $F=G=0$. Likewise the only
solution of the equation $e^{it\xi^3}\overline{F(\xi )} =\overline{G(\xi )}$ with
$F,G\in H^2$ is $F=G=0.$\par
\indent Next prove that all non-zero vectors in $W_tL^2\cap L^2({\bf R}_+)$ are cyclic. 
Suppose that $G\neq 0$ is a non-cyclic vector in $H^2\cap \overline{ e^{it\xi^3}H^2}$ so
that $\overline{G(\xi )}=e^{it\xi^3}F(\xi )$ for some $F\in H^2$, and where $G\perp
uH^2$ for some inner function $u$. We have $u\overline{G}\in H^2$; so we introduce
inner functions $v$ and $w$, and an outer function $\theta$, such that $u\overline
{G}=v\theta$ and $F=w\theta.$ Then, as in [9, Theorem 3.1.1],
$$e^{it\xi^3}={{\overline{G}}\over {F}}={{v}\over {uw}}\eqno(5.14)$$
\noindent is a quotient of inner functions and hence is of finite Nevanlinna
type, but the corresponding logarithmic integral (5.12) diverges, and we have a
contradiction. (The author conjectures that $W_tL^2\cap L^2({\bf R}_+)=0$ so that 
$W_tL^2$ and $L^2({\bf R}_+)$ are in general position, since any non-zero elements in the
intersecting subspaces would satisfy some implausible equations.)\par
\indent (ii) The Hankel operator $\Gamma_{0,t}=P_+J_tP_+=P_+e^{itD^3}JP_+$  has
kernel
$${{1}\over{(3t)^{1/3}}} {\bf I}_{(0,\infty )}(x){\hbox{Ai}}
\Bigl( {{x+y}\over 
{(3t)^{1/3}}}\Bigr) {\bf I}_{(0,\infty )}(y),\eqno(5.15)$$
\noindent which is of Hilbert--Schmidt type; see
[{23}, page 46]. Indeed,
the integral
$$\int_0^\infty\!\!\!\int_0^\infty {{1}\over {(3t)^{2/3}}}{\hbox{Ai}}
\Bigl( {{x+y}\over 
{(3t)^{1/3}}}\Bigr)^2 dxdy$$
\noindent may be transformed by the substitution $u=x+y$ to the convergent integral
$${{1}\over {(3t)^{2/3}}}\int_0^\infty u{\hbox{Ai}}\Bigl( {{u}\over
{(3t)^{1/3}}}\Bigr)^2\, du<\infty;\eqno(5.16)$$
\noindent here we use the bounds from [11, page 43] 
$${\hbox{Ai}}(x)={{1}\over {2\sqrt\pi x^{1/4}}}\bigl(
1+O(x^{-3/2})\bigr)\exp\Bigl(-{{2}\over{3}}x^{3/2}\Bigr)\qquad
(x\rightarrow\infty ).\eqno(5.17)$$
\indent Hence the off-diagonal 
operators $P_-W_tP_+=P_-J_t(P_+J_tP_+)$ and\par
\noindent 
$P_+W_tP_-=(P_+J_tP_+)J_tP_-$
are Hilbert--Schmidt. For the bottom-right entry, we have a stronger
conclusion, namely that $P_+W_tP_+=(P_+J_tP_+)(P_+J_tP_+)$ is trace
class.\par
\indent (iii) When we replace $t\geq 0$ by $t\leq 0$, we need to switch the
roles of $P_+$ and $P_-$ in the previous discussion and we deduce that 
$P_-W_tP_+$ and $P_+W_tP_-$ are Hilbert--Schmidt, while $P_-W_tP_-$ is of trace class. \par
\rightline {$\square$}\par
\vskip.05in
\noindent {\bf Theorem 5.4.} {\sl (i) The $C_0$ unitary groups 
$S_s$ and $U_t=e^{-it(D-x^2)}$ satisfy the Weyl relations 
$S_sU_t=e^{ist}U_tS_s$ for $s,t\in {\bf R}.$\par 
\indent (ii) For $\alpha\geq 0$ and real $\delta$, the subspace $e^{ix^3/3-i\alpha
x^2+i\delta x}H^2$ is simply invariant for $S_s$ $(s\geq 0)$ and $U_t$
$(t\geq 0)$. Conversely, if ${\cal T}$ is a non-trivial simply invariant subspace 
for $S_s$ $(s\geq 0)$ and 
for $U_s$
$(s\geq 0)$, then ${\cal T}=e^{ix^3/3-i\alpha
x^2+i\delta x}H^2$ for some $\alpha\geq 0$ and real $\delta .$\/}\par
\vskip.05in
\noindent {\bf Proof.} (i) One can prove directly that the operators $U_t$ defined by 
$$U_tf(x)=e^{i(x^2t-xt^2+t^3/3)}f(x-t)\eqno(5.18)$$
\noindent define a $C_0$ unitary group on $L^2({\bf R}).$
This formula for the $U_t$ of (5.18) was obtained by the method of characteristics.
Indeed, when $f$ is differentiable, the function $g(x,t)=e^{i(x^3-(x-t)^3)/3}f(x-t)$
satisfies
$$\eqalignno{{{\partial g}\over {\partial t}}+{{\partial g}\over {\partial
x}}&=ix^2g,\cr
g(x,0)&=f(x);&(5.19)\cr}$$
\noindent and so $U_tf(x)=e^{-it(D-x^2)}f(x)=g(x,t)$ gives the unique solution of 
the initial value problem (5.19).\par
\indent Let $V$ be the unitary operator $V:f(x)\mapsto e^{ix^3/3}f(x)$ on
$L^2({\bf R})$, then clearly $S_s=V^*S_sV.$ The generator of the unitary group 
$V^*U_tV$ equals 
$$-iV^*(D-x^2)V=-ie^{-ix^3/3}(D-x^2)e^{ix^3/3}=-{{\partial}\over
{\partial x}}=-iD;\eqno(5.20)$$
so by the uniqueness of groups with given generator we have $V^*U_sV=e^{-isD}=\tau_s$
and hence $U_s=Ve^{-isD}V^*=V\tau_sV^*$. By conjugating the Weyl relations 
$\tau_sS_t=e^{-ist}S_t\tau_s$ for $(s,t\in {\bf R})$ by $V$, we can deduce (5.19).\par
\indent (ii) Clearly any ${\cal T}=e^{ix^3/3-i\alpha x^2+i\delta x}H^2$ is
simply invariant under $S_s$
$(s\geq 0)$, and we can use the preceding calculations to show that ${\cal T}$ is also
invariant for $U_s$ $(s\geq 0)$. Indeed, for $g\in {\cal T}$ we can take 
$f\in H^2$ such that $g(x)=e^{ix^3/3-i\alpha x^2+i\delta x}f(x)$ and we have
$$\eqalignno{U_sg&=U_s(e^{ix^3/3-i\alpha x^2+i\delta x}f)\cr
&=V\tau_sV^*V \{e^{-i\alpha x^2+i\delta x}f\}\cr
& =V\tau_s\{ e^{-i\alpha x^2+i\delta x}f\}\cr
&=e^{2i\alpha sx-i\alpha s^2-i\delta s} e^{ix^3/3-i\alpha x^2+i\delta
x}f(x-s)&(5.21)\cr}$$
\noindent where $f(x-s)$ is an $H^2$ function, so $U_sg$ belongs to 
the subspace $S_{2\alpha s}{\cal T}$ of ${\cal T}$. This proves the forward implication. \par
\indent To prove the converse, we take any ${\cal T}$ that is simply invariant as in the
Theorem, and observe that $V^*{\cal T}$ is simply invariant under $S_s$ $(s\geq 0)$ since
$V^*$ commutes with $S_s$, and $V^*{\cal T}$ is also invariant under $\tau_s$
 $(s\geq 0)$ since 
$\tau_sV^*{\cal T}=V^*U_s{\cal T}\subseteq V^*{\cal T}.$ By the Katavolos--Power 
Theorem [{15}], 
there exist $\alpha>0$ and a real $\delta$ such that $V^*{\cal T}=e^{-i\alpha x^2+i\delta x}H^2,$ and
hence ${\cal T}$ has the required form.\par
\rightline {$\square$}\par
\noindent {\bf Corollary 5.5.} {\sl The unitary groups $\tau_{-s}=e^{isD}$ and $\hat 
U_t=e^{it(D^2+x)}$ form a Weyl pair on $L^2$, such that the lattice of subspaces $\{ 
\tau_{-\delta }\hat U_{-\alpha }W_{1/3}L^2: \delta \in {\bf R}, \alpha >0\}$ gives
the set of simply
invariant closed linear subspaces for $(\tau_{-s})$ $(s\geq 0)$ and $(\hat U_t)$ $(t\geq 0).$}\par
\vskip.05in
\noindent {\bf Proof.} The strongly continuous unitary semigroups 
 $\tau_{-s}={\cal F}^*S_s{\cal F}$ $(s\geq 0)$ and 
${\cal F}^*U_t{\cal F}$
and ${\cal F}^*S_t{\cal F}$ $(t\geq 0)$ have generators $iD$ and $i(D^2+x)$
respectively, so $(\tau_{-s}, \hat U_t)$ forms a Weyl pair by 
Theorem 5.4(ii) and Proposition 5.2(i). Under the Fourier transform, we have
$$\eqalignno{{\cal F}\tau_{-\delta}\tau_{-\alpha^2}\hat U_{-\alpha}W_{1/3}L^2&=S_\delta
S_{\alpha^2}U_{-\alpha}VH^2\cr
&=\{ e^{i\delta x-i\alpha x^2 +ix^3/3}f(x): f\in H^2\};&(5.22)\cr}$$
\noindent hence the subspaces correspond as stated in the Corollary.\par
\vskip.1in
\noindent {\bf 6. Mathieu functions and periodic potentials}\par
\indent One can construct Tracy--Widom operators over the circle by means of the
differential equations of section 2. Kernels of the form (6.6) below arise in the
theory of random unitary matrices with Haar probability measure as in [22, p. 195]. The
purpose of this section is to introduce examples beyond the list in [30].\par
\indent Let $S$ be the $2\times 2$ fundamental solution matrix of Hill's equation with
smooth $\pi$-periodic potential $q$, so that 
$${{d}\over {dx}} S=\left[\matrix{0& 1\cr -(\lambda +q(x))& 0\cr}\right] S, \qquad
S(0)=\left[\matrix{1&0\cr 0&1\cr}\right];\eqno(6.1)$$
\noindent then $\det S=1$, and $\Delta (\lambda )={\hbox {trace }}\, S(\pi )$ defines
the discriminant. When $\lambda$ is real, evidently $\Delta (\lambda
)^2\geq 4$ if and only if the eigenvalues of $S(\pi)$ are real, and $\Delta (\lambda )^2=4$ occurs
if and only if $S$ is periodic with period $\pi$ or $2\pi$. The periodic spectrum 
$$\Lambda =\{\lambda_0< \lambda_1\leq \lambda_2<\lambda_3\leq \lambda_4<\dots
<\lambda_n\nearrow \infty \}\eqno(6.2)$$
\noindent of Hill's equation consists of those real $\lambda$ such that 
$$y''+(\lambda +q)y=0\eqno(6.3)$$
\noindent has a non-trivial $\pi$ or $2\pi$ periodic solution as in [18, p.11]. The discriminant
satisfies 
$$4-\Delta^2 (\lambda )=4(\lambda -\lambda_0)\prod_{j=1}^\infty
{{(\lambda_{2j-1}-\lambda )(\lambda_{2j}-\lambda)}\over {j^4}}.\eqno(6.4)$$ 
\vskip.05in
\noindent {\bf Theorem 6.1} {\sl For each real $\alpha$ there exists an infinite
sequence of $\lambda_n$ such that Hill's equation (6.3) with potential $\alpha \cos
2x$ has a non-trivial $2\pi$-periodic and real solution $A$.
For such $A$, let $W$ be the kernel 
$$W(x,y) ={{A(x)A'(y)-A'(x)A(y)}\over {\sin (x-y)}},\eqno(6.5)$$
\noindent which continuously differentiable and doubly periodic with period
$2\pi $.\par
\indent (i) Then $W$ defines a self-adjoint and Hilbert--Schmidt operator on
$L^2[0, 2\pi ]$;\par 
\indent (ii) the eigenfunction corresponding to each non-zero simple eigenvalue of $W$
is a $2\pi$-periodic solution of (6.3).}\par
\vskip.05in
\noindent {\bf Proof.} (i) When $\alpha =0$ and $\lambda =n^2$ with $n=1, 2, \dots $, we can take
$A(x)=\sin nx$, and recover the kernel
$$W(x,y)={{n\sin n(x-y)}\over {\sin (x-y)}}\eqno(6.6)$$
\noindent as in the circular ensemble.\par
\indent When $\alpha \neq 0$, there exists by Hochstadt's theorem [18, page 40] an increasing sequence $(\lambda_n')$
which satisfies the estimates  
$$\lambda_{2n-1}'= (2n-1)^2 +{{\alpha^2}\over {32n^2}} +o(n^{-2})\qquad
(n\rightarrow\infty )\eqno(6.7)$$
$$0<\lambda_{2n}'- \lambda_{2n-1}'=o(n^{-2}),$$
\noindent and such that (6.3) has a non-trivial solution $A$. This is Mathieu's function
of the first kind. As in section 2, we can calculate
$$\Bigl({{\partial }\over {\partial x}} +{{\partial }\over {\partial y}}\Bigr) W(x,y)=
-2\alpha (\sin x \cos y+\cos x\sin y) A(x)A(y),\eqno(6.8)$$
\noindent hence
$$ W(x,y)=-{\alpha}\int_0^{x+y}\sin \theta \,A\Bigl(
{{1}\over{2}}(\theta +x-y)\Bigr) A\Bigl(
{{1}\over{2}}(\theta -x+y)\Bigr)\,
d\theta +g(x-y),\eqno(6.9)$$
\noindent where, by letting $x=-y$, one finds 
$$g(x)={{A(x/2)A'(-x/2)-A'(x/2)A(-x/2)}\over {\sin x}}.$$
\noindent Evidently $W$ is a real, symmetrical and continuous kernel, and hence
determines a self-adjoint and Hilbert--Schmidt operator on $L^2[0,2\pi ].$\par 
\indent (ii) By differentiating (6.8) and recalling the definition of $W$, one can easily deduce that 
$$\Bigl( {{\partial^2}\over {\partial x^2}}- {{\partial^2}\over {\partial y^2}}\Bigr)
W(x,y)=\alpha (\cos 2x-\cos 2y) W(x,y).\eqno(6.10)$$
\indent For $\nu\neq 0$, any non-zero solution $f\in L^2[0,2\pi ]$ of the integral equation 
$$\nu f(x)=\int_0^{2\pi} W(x,y)f(y)\, dy\eqno(6.11)$$ 
\noindent extends to define a twice continuously differentiable and $2\pi$-periodic function on ${\bf R}$. Now $g(x) 
=f''(x)+\alpha \cos 2x f(x)$ also gives a $2\pi$ periodic and continuous
solution of (6.11); this follows from (6.10) by an
integration-by-parts argument. By simplicity of the eigenvalue, we deduce
that $g$ is a constant multiple of $f$, and hence that $f$ is a $2\pi$ periodic
solution of Mathieu's equation.\par

\rightline {$\square$}\par
\vskip.05in
\indent Conversely, let $M_\Lambda$ be the
space of potentials $q$ that have periodic spectrum equal to a given $\Lambda$. McKean,
van Moerbeke and Trubowitz [19, 20] have shown that $M_\Lambda$ can be considered as a torus 
$$M_\lambda =\Bigl\{ {{1}\over{2}}\Bigl(\Delta (x_j)+\sqrt{\Delta
(x_j)^2-4}\Bigr)_{j=1}^\infty : \lambda_{2j-1}\leq x_j\leq \lambda_{2j};
j=1,2,\dots \Bigr\}\eqno(6.12)$$ 
\noindent over the product over the intervals of instability 
$(\lambda_{2j-1},\lambda_{2j})$ where $\Delta (\lambda
)^2<4$ and that $M_\Lambda$ is associated
with the Jacobi manifold over the Riemann surface of $\sqrt{\Delta^2(\lambda )-4}$. Hence 
$M_\lambda$ can have dimension $n=0, 1, \dots, \infty$, equal to the number of
simple zeros 
$\Delta (\lambda )^2-4$. The periodic spectrum $\Lambda$ is preserved by
Hamiltonian flows; in particular, there is a $2\pi$ periodic Korteweg--de Vries flow on 
$M_\Lambda$ associated with
$${{\partial q}\over {\partial t}}=3q{{\partial q}\over {\partial x}}
-{{1}\over {2}}{{\partial^3 q}\over {\partial x^3}}.\eqno(6.13)$$

\noindent By Theorem 6.1, the potential $\alpha\cos 2x$ gives an infinite-dimensional $M_\Lambda$ on which there are solutions to $KdV$ that are $2\pi$ periodic
in the space variable and almost periodic in time [4, Appendix]. This makes an 
interesting contrast with the concentric KdV equation $u_t+u/(2t)-6uu_x+u_{xxx}=0$
which is used via $\Gamma_{(x)}^2$ in Theorem 2.3 to linearize
${\hbox{P}}_{{\hbox{\smalletters II}}}$ for the soft edge
ensemble [10, page 173].\par

\vskip.05in
\noindent {\bf Acknowledgements}\par
\noindent I am grateful to Stephen Power and Sergei Treil for helpful conversations. This work was partially supported by EU 
Training Network Grants HPRN-CT-2000-00116 `Classical Analysis and Operators' and 
MRTN-CT-2004-511953 `Phenomena in High
Dimensions'.\par
\vskip.1in
\noindent {\bf References}\par
\noindent [1] M.J. Ablowitz and P.A. Clarkson, {Solitons, Nonlinear Evolution Equations 
and Inverse Scattering,} Cambridge University Press, Cambridge, 1991.\par
\noindent [2] D. Alpay, {The Schur algorithm, reproducing kernel Hilbert spaces,
and system theory,} American Mathematical Society/Soci\'et\'e Math\'ematique de France,
2001.\par
\noindent [3] G. Aubrun, {A sharp small deviation inequality for the largest
eigenvalue of a random matrix,}
S\'eminaire de Probabilit\'es XXXVIII 320--337, Lecture Notes in Math. {1857},
Springer, Berlin, 2005.\par
\noindent [4] J. Bourgain, Fourier transform restriction phenomena for certain lattice
subsets and applications to nonlinear evolution equations II: the KdV equation, {
Geom. Funct. Anal.} {3} (1993), 209--262.\par
\noindent [5] A. Borodin, {Biorthogonal ensembles,}  
{Nuclear Phys. B}  {536}  (1999), 704--732.\par
\noindent [6] L. De Branges, {Hilbert Spaces of Entire Functions}, Prentice--Hall,
Englewood Cliffs, 1968.\par
\noindent [7] J.--F. Burnol, Sur les `espaces de Sonine' associ\'es par de Branges \`a
la transformation de Fourier, {C.R. Math. Acad. Sci. Paris} {335} (2002),
689--692.\par  
\noindent [8] P.A. Deift, A.R. Its and X. Zhou, 
A Riemann--Hilbert approach to asymptotic problems arising in the theory of 
random matrix models, and also in the theory of integrable statistical mechanics,  
{Ann. of Math.} (2)  {146}  (1997), 149--235.\par
\noindent [9] R.G. Douglas, H.S. Shapiro and A.L. Shields, Cyclic vectors and invariant
subspaces for the backward shift operator, {Ann. Inst. Fourier (Grenoble)} {20} 
(1970), 37--76.\par
\noindent [10] P.G. Drazin and R.S. Johnson, {Solitons: an Introduction}, Cambridge
University Press, Cambridge, 1989.\par
\noindent [11] A. Erd\'elyi, {Asymptotic Expansions}, Dover, California, 1955.\par
\noindent [12] P.J. Forrester, The spectrum edge of random matrix ensembles, {Nuclear
Physics B} {402} (1993), 709-728.\par
\noindent [13] P.J. Forrester and E.M. Peter J.; Rains, 
Correlations for superpositions and decimations of Laguerre and Jacobi orthogonal 
matrix ensembles with a parameter,  Probab. Theory Related Fields  130  (2004), 
518--576. \par
\noindent [14] E.L. Ince, {Ordinary differential equations},  
Dover Publications, London, 1956.\par
\noindent [15] A. Katavolos and S.C. Power, {The Fourier binest
algebra,} {Math. Proc. Cambridge Philos. Soc.} {122} (1997), 
525--539.\par
\noindent [16] A. Katavolos and S.C. Power, {Translation and dilation invariant
subspaces of $L^2({\bf R})$,} {J. reine angew. Math.} {552} (2002), 
101--129.\par
\noindent [17] P. Koosis, {Introduction to $H_p$ Spaces},
Cambridge University Press, Cambridge, 1980.\par
\noindent [18] W. Magnus and S. Winkler, {Hill's Equation}, Dover, New York,
1966.\par
\noindent [19] H.P. McKean and P. van Moerbecke, The spectrum of Hill's equation, 
{Invent. Math.} {30} (1975), 217--274.\par
\noindent [20] H.P. McKean and E. Trubowitz, Hill's operator and hyperelliptic
function theory in the presence of infinitely many branch points, {Comm. Pure Appl. 
Math.} {29} (1976), 143--226.\par  
\noindent [21] A.N. Megretski\u\i, V.V. Peller and S.R. Treil, 
The inverse spectral problem for self-adjoint Hankel operators, {Acta Math.} { 
174} (1995), 241--309.\par
\noindent [22] M.L. Mehta, {Random Matrices,} 2nd ed.,
Academic Press, San Diego, 1991.\par
\noindent [23] V. Peller, {Hankel Operators and Their
Applications}, Springer, New York, 2003.\par
\noindent [24] J. Rovnyak and V. Rovnyak, {Sonine spaces of entire functions}, J. Math.
Anal. Appl. 27 (1969), 68--100.\par
\noindent [25] I. N. Sneddon, {The use of Integral Transforms,} McGraw--Hill, New Delhi,
1974.\par
\noindent [26] N. Sonine, Recherches sur les fonctions cylindriques et le
d\'eveloppement des fonctions continue en s\'eries, Math. Annal. 16 (1880), 1--80.\par 
\noindent [27] G. Szeg\"o, {Orthogonal Polynomials,}
American Mathematical Society, New York, 1959.\par
\noindent [28] C.A. Tracy and H. Widom, {Level spacing distribution and the Airy
kernel,} {Comm. Math. Phys.} {159} (1994), 
151--174.\par
\noindent [29] C.A. Tracy and H. Widom, {Level spacing distribution and the Bessel
kernel,} {Comm. Math. Phys.} {161} (1994), 
289--309.\par
\noindent [30] C.A. Tracy and H. Widom, {Fredholm determinants, differential
equations and matrix models}, {Comm. Math. Phys.} {163} (1994), 33--72.\par
\noindent [31] C.A. Tracy and H. Widom, {A system of differential equations for the
Airy process}, {Electron. Comm. Probab.} {8} (2003), 93--98.\par
\noindent [32] J. von Neumann, {Die Eindeutigkeit der Schr\"odinger
Operatoren,} {Math. Ann.} {104} (1931), 570--578.\par
\vfill
\eject
\end